\documentstyle[12pt,leqno]{amsart}

\font\chq = cmr7
 0

\numberwithin{equation}{section}
\theoremstyle{plain}
	\begingroup
	\newtheorem{thm}[equation]{Theorem}
	\newtheorem{lem}[equation]{Lemma}
	\newtheorem{prop}[equation]{Proposition}
	\newtheorem{cor}[equation]{Corollary}
	
\endgroup
\theoremstyle{remark}
	\newtheorem{rem}[equation]{Remark}
\theoremstyle{definition}
	\begingroup
	\newtheorem{defn}[equation]{Definition}
	\newtheorem{exmp}[equation]{Example}
	
	\newtheorem{algo}[equation]{}
	\newtheorem{tactic}[equation]{Tactic}
\endgroup

\newcommand{\so}{\Rightarrow}

\newsymbol\twoheadrightarrow 1310
\newcommand{\sobre}{\twoheadrightarrow}             
\newcommand{\inclu}{\hookrightarrow}
\newcommand{\ganda}{\rightharpoonup}                

\newcommand{\End}{\operatorname{End}}
\newcommand{\Aut}{\operatorname{Aut}}

\newcommand{\GL}{\operatorname{GL}}

\newcommand{\Ad}{\operatorname{Ad}}

\newcommand{\id}{\operatorname{id}}

\newcommand{\Ind}{\operatorname{Ind}}

\newcommand{\car}{\operatorname{char}}

\newcommand{\gr}{\operatorname{gr}}
\newcommand{\Ss}{{\mathcal S}}

\newcommand{\ZZ}{{\Bbb Z}}
\newcommand{\NN}{{\Bbb N}}
\renewcommand{\SS}{{\Bbb S}}

\newcommand{\DD}{{\Bbb D}}
\newcommand{\HH}{{\Bbb H}}

\renewcommand{\k}{{\mathbf k}}

\newcommand{\toba}[1]{{\frak t}(#1)}                
\newcommand{\tobag}[2]{{\frak t}^{#1}(#2)}      

\newcommand{\cate}{{\mathcal C}}                    
\newcommand{\cateb}{\overline{{\mathcal C}}}        
\newcommand{\yd}{{}^H_H{\mathcal YD}}               
\newcommand{\yds}[1]{{}^{#1}_{#1}{\mathcal YD}}     
\newcommand{\ydskg}{{}^{\k\Gamma}_{\k\Gamma}{\mathcal YD}} 
\newcommand{\coind}[2]{#1^{\hbox{\chq co}#2}}       
\newcommand{\dad}[1]{#1^*}                          
\newcommand{\agg}[2]{#1(#2)}                        

\newcommand{\prim}[1]{{\mathcal P}(#1)}         
\newcommand{\primg}[2]{{\mathcal P}_{#2,1}(#1)} 
\newcommand{\primgh}[3]{{\mathcal P}_{#2,#3}(#1)}   
\newcommand{\unid}[1]{#1^{\times}}              
\newcommand{\ii}{{\underline i}}                
\newcommand{\jj}{{\underline j}}                
\newcommand{\kk}{{\underline k}}                
\newcommand{\alai}[2]{{}^{#2}#1}                
\newcommand{\com}[2]{#1_{(#2)}}                 
\newcommand{\con}[2]{#1_{(-#2)}}                
\newcommand{\ctr}[3]{{#1 \choose #2}_{\!\!\!#3}\,}  
\newcommand{\fmu}{\phantom{{}^{-1}}}		
\newcommand{\fsm}{\phantom{-}}			
\newcommand{\mdpd}[4]{\left(\begin{array}{rr}%
	#1 & #2 \\ #3 & #4\end{array}\right)}   

\newcommand{\noi}{\noindent}

\newcommand{\Gm}{\Gamma}
\newcommand{\Dt}{\Delta}
\newcommand{\dt}{\delta}
\newcommand{\sgm}{\sigma}
\newcommand{\te}{\theta}
\newcommand{\al}{\alpha}

\newcommand{\lmb}{\lambda}
\newcommand{\ptl}{\partial}

\newcommand{\rmu}{{\bold i}}
\renewcommand{\toba}[1]{{\frak B}(#1)}
\renewcommand{\tobag}[2]{{\frak B}^{#1}(#2)}

\newcommand{\ydm}{YD-module}
\newcommand{\ceun}{\{0,1\}}
\newcommand{\primgp}[1]{{\mathcal P}'_{#1,1}} 
\newcommand{\primghp}[2]{{\mathcal P}'_{#1,#2}} 
\numberwithin{equation}{subsection}
\renewcommand{\theequation}{\arabic{section}.\arabic{subsection}.\arabic{equation}}

\begin{document}

\renewcommand{\baselinestretch}{1.2}
\renewcommand{\thefootnote}{}
\thispagestyle{empty}
\vspace*{2in}
\begin{center}
{\Large \bf \sc Pointed Hopf algebras of dimension 32} \\
\vspace{2em}
{\large Mat\'\i as Gra\~na} \\
\vspace{1em}
Depto de Matem\'atica, Pab. I, Ciudad Universitaria \\
(1428) Buenos Aires, Argentina \\
e-mail: {\tt matiasg@@dm.uba.ar}

\footnote{This work was partially supported by CONICET,
	CONICOR, Secyt (UNC), Secyt (UBA)}

{Abstract}
\end{center}
\small
We give a complete classification of the $32$-dimensional pointed
Hopf algebras over an algebraically closed field $\k$ with
$\car\k\neq 2$.
It turns out that there are infinite families of isomorphism
classes of pointed Hopf algebras of dimension $32$.
In \cite{as1}, \cite{bdg} and \cite{ge} are given families
of counterexamples for the tenth Kaplansky conjecture.
Up to now, $32$ is the lowest dimension where Kaplansky conjecture fails.

\normalsize
\section{Introduction}

In recent years there has been certain activity in classification problems
of pointed Hopf algebras.
In \cite{as1}, \cite{bdg}, \cite{ge} independently the authors
find infinite families of pointed Hopf algebras of dimension $p^4$
($p$ an odd prime number), solving thus by the negative the tenth
Kaplansky conjecture. In \cite{dso} there is a classification of
pointed Hopf algebras of dimension $p^n$ with coradical of dimension
$p^{n-1}$.

In \cite{as1}, \cite{as2}, \cite{as3} there are classifications
of pointed Hopf algebras of dimension $p^3$, $p^4$ and a large number of
strong results towards the classification of pointed Hopf algebras with
commutative coradical. Using the techniques of these articles,
namely the ``lifting procedure" described below,
we classify here the pointed Hopf algebras of dimension $32$ over
an algebraically closed field of characteristic $\neq 2$. Related results
can be found in \cite{cdr} (the classification of the pointed Hopf algebras
of dimension $16$), in \cite{be} (classification of Ore extension Hopf algebras
with coradical $\k(C_2\times C_2)$ and $\k C_4$), and in \cite{gp5}
(the description of the coradically graded pointed Hopf algebras of
dimension $p^5$).

One of the key results is Corollary \ref{pr:cto}, where we compute
the Nichols algebras $\toba V$ (see \ref{al:prl} for the definition
of Nichols algebras), $V$ the Yetter--Drinfeld modules of dimension
$2$ over the groups of order $2^j$, $1\le j\le 4$.

\bigskip
\noindent Our main result is
\begin{thm}
Let $A$ be a pointed Hopf algebra of dimension $32$ over the algebraically
closed field $\k$, $\car\k\neq 2$. Let $\Gm=G(A)$ be the group of
group-likes of $A$.

\noi If $\Gm=C_2$ then $A$ is the algebra of section \ref{sz2}.

\noi If $\Gm=C_2\times C_2$ then $A$ is one of the $6$ algebras of \ref{sz2z2}.

\noi If $\Gm=C_4$ then $A$ is one of the $12$ algebras of \ref{sz4}.

\noi If $\Gm=C_2\times C_2\times C_2$ then $A$ is one of the $6$ algebras of \ref{sz2z2z2}.

\noi If $\Gm=C_2\times C_4$ then $A$ is in the infinite family of algebras of \ref{sz2z4}.

\noi If $\Gm=C_8$ then $A$ is in the infinite family of algebras of \ref{sz8}.

\noi If $\Gm=\DD_4$ the dihedral group, then $A$ is one of the $7$ algebras of \ref{sd4}.

\noi If $\Gm=\HH$ the quaternionic group, then $A$ is one of the $3$ algebras of \ref{sh}.

\noi If $\Gm=C_2\times C_2\times C_2\times C_2$ then $A$ is the algebra of \ref{sz2z2z2z2}.

\noi If $\Gm=C_2\times C_2\times C_4$ then $A$ is one of the $7$ algebras of \ref{sz2z2z4}.

\noi If $\Gm=C_4\times C_4$ then $A$ is one of the $4$ algebras of \ref{sz4z4}.

\noi If $\Gm=C_2\times C_8$ then $A$ is one of the $14$ algebras of \ref{sz2z8}.

\noi If $\Gm=C_{16}$ then $A$ is one of the $7$ algebras of \ref{sz16}.

\noi If $\Gm$ is one of the non abelian groups of order $16$ then $A$ is one of
	the $13$ algebras of \ref{sgna}.

\noi If $\Gm$ has order $32$ then $A$ is a group algebra.
\end{thm}

\subsection{Some basic definitions and conventions}\label{ss:bdac}
Throughout the article $\k$ shall be an algebraically
closed field of characteristic $\neq 2$.
The degrees for graded algebras shall be denoted by parentheses, i.e.
$\agg Hi$ is the homogeneous component of degree $i$ of the graded
algebra $H$.
The symbols $\rmu,\xi$ shall denote fixed roots of unity of orders
$4,8$ respectively, such that $\rmu=\xi^2$.

Let $H$ be a Hopf algebra. We denote by $\yd$ the category of (left-left)
Yetter--Drinfeld modules over $H$, i.e. $M$ is an object of $\yd$
if it is a left
$H$-module (with structure $h\otimes m\mapsto h\ganda m$),
a left $H$-comodule (with structure $\dt(m)=\con m1\otimes \com m0$) and
such that
$$\dt(h\ganda m)=\com h1\con m1\Ss(\com h3)\otimes \com h2\com m0.$$
We denote by $c$ the braiding in $\yd$, which is given by
$$c=c_{M,N}:M\otimes N\to N\otimes M,\quad
	m\otimes n\mapsto \con m1\ganda n\otimes \com m0.$$
When $H$ is clear from the context, we shall call a module
in $\yd$ a \ydm.
Unless explicitly stated, the letter $V$ shall denote a \ydm.  

For $\Gm$ a group and $g\in\Gm$ we denote by $\Gm_g$ the commutator
subgroup
$$\Gm_g=\{h\in\Gm\ |\ hg=gh\}.$$
It is not difficult to prove (see for instance \cite{ag}) that if
$H=\k\Gm$ and $|\Gm|\neq 0$ in $\k$ then $\yd$ is a semisimple category
and the irreducible objects are given as follows: let $g\in\Gm$ and
$\rho:\Gm_g\to\Aut(W)\in\widehat{\Gm_g}$ an irreducible representation.
Then the induced representation
$$M(g,\rho)=\Ind_{\Gm_g}^{\Gm}W=\k\Gm\otimes_{\k\Gm_g}W$$
with the usual module structure and the comodule structure given by
$$\dt(h\otimes w)=hgh^{-1}\otimes (h\otimes w)\in\k\Gm\otimes M(g,\rho)$$
is an irreducible YD-module. Furthermore, let $\{g_1,\ldots,g_n\}$ be
a subset of $\Gm$ consisting of one element for each conjugacy class,
and let $\{\rho_i^1,\ldots,\rho_i^{m_i}\}$ be the irreducible representations
of $\Gm_{g_i}$. Then the objects $M(g_i,\rho_i^j)$ ($1\le i\le n$,
$1\le j\le m_i$) gives a full collection of
non-isomorphic irreducible objects of $\ydskg$.

We extend the definition of $M(g,\rho)$ to the case in which
$\rho$ is not irreducible. Observe that
$\dim M(g,\rho)=\deg\rho\times[\Gm:\Gm_g]$. Hence, if
$\Gm$ is abelian, the irreducible modules of $\ydskg$ are
$1$-dimensional. In this case, let $V=\oplus_iM(h_i,\chi_i)$ where
$h_i\in\Gm$ and $\chi_i\in\hat\Gm$, and denote by $x_i$ a generator
of $M(h_i,\chi_i)$. By the definition of the braiding,
\begin{equation}\label{eq:bca}
c(x_i\otimes x_j)=\chi_j(h_i)x_j\otimes x_i.
\end{equation}

Let $q\in\k$. For $n\ge m\in\NN$, we use the standard notation
$$(n)_q=\sum_{i=0}^{n-1}q^i,\quad
(n)^{!}_q=\prod_{i=1}^n(i)_q,\quad
\ctr nmq=\frac{(n)^{!}_q}{(m)^{!}_q(n-m)^{!}_q}.$$

\medskip
\noindent Let $H$ be a graded Hopf algebra, $H=\oplus_i\agg Hi$.
Let $\pi:H\to\agg H0$ and $\iota:\agg H0\to H$ be respectively the
canonical projection and the canonical inclusion. Then
$$R=\coind H{\pi}=\coind H{\agg H0}
	=\{x\in H\ |\ (\id\otimes\pi)\Dt(x)=x\otimes 1\}$$
is a braided Hopf algebra in $\yds{\agg H0}$. Furthermore, it is
graded: $R=\oplus_i\agg Ri$, with $\agg Ri=R\cap\agg Hi$. The algebra
$H$ can be reconstructed as a bosonization, or biproduct
(see e.g. \cite{ag})
$$H=R\#\agg H0.$$

A braided Hopf algebra $R$ in $\yds K$ ($K$ a Hopf algebra) is called
{\it Nichols algebra} if
\begin{align}
&R=\oplus_i\agg Ri\mbox{ is a graded braided Hopf algebra},\label{toba1} \\
&\agg R0=\k 1,\label{toba2} \\
&\prim R=\agg R1,\label{toba3} \\
&\agg R1\mbox{ generates }R,\label{toba4}
\end{align}
where $\prim R$ stands for the space of primitive elements of $R$.
It is not known whether in characteristic $0$ and for $R$ finite
dimensional the condition \eqref{toba4}, resp. \eqref{toba3}, is superfluous,
though there are some results in this direction. For instance,
in \cite[Thm. 3.2]{as1} is proved that if $\dim\agg R1=1$ ($R$ finite
dimensional and $\car\k=0$) then
(\eqref{toba1}--\eqref{toba3})$\so$\eqref{toba4}.
A Nichols algebra $R$ is uniquely determined by its primitive elements
$\prim R$. We use hence the notation $R=\toba V$ for the Nichols algebra
$R$ with $\agg R1=V$, whenever
$V\in\yds K$. We specify sometimes the graduation of $\toba V$ by
$$\toba V=\oplus_i\tobag iV.$$

The construction $V\mapsto\toba V$ can be made with the so called
{\it quantum-shuffles} (see for instance \cite{ag}), though this
construction does not depend on the braided category $\yds K$, but
only on the pair $(V,c)$. The quantum-shuffles construction applied to
any pair $(V,c)$ ($c\in\End(V\otimes V)$ a solution of the braid equation),
the result has been called {\it Quantum Shuffle Algebra} in \cite{ros}.
Nichols algebras have some particularly useful properties which
distinguish them from this general concept of Quantum Shuffle Algebras,
e.g. the kind of Poincar\'e duality stated in \ref{pr:dualpuan}.

\medskip
We shall use the following notation for finite abelian groups: let
$\Gm=C_{n_1}\times\cdots\times C_{n_t}$. We denote by $g_j$
a generator of the $j$-th factor, i.e. $\Gm$ can be presented by
generators $g_j$ and relations $g_j^{n_j}=[g_j,g_k]=1$. Furthermore, let
$\al_j=-1,\rmu,\xi,\sqrt{\xi}$ depending whether $n_j=2,4,8,16$ (recall
that $\xi^2=\rmu,\ \rmu^2=-1$). We denote by
$\hat g_i$ the character of $\Gm$ given by
$\hat g_i(g_j)=\al_i^{\dt_{ij}}$. Then
$\hat\Gm\simeq\Gm$ with generators $\hat g_i,\ i=1,\ldots,t$.

If $x_1,\ldots,x_r$ are vectors in a space $V$, we denote by
$(x_1,\ldots,x_r)$ the linear span of them. If $g_1,\ldots,g_r$
are elements in $\Gm$, we denote by $(g_1,\ldots,g_r)$ the subgroup
generated by them.

If $V\in\ydskg$ and $\Gm$ is abelian then
$V=\oplus_{i=1}^{\te}M(h_i,\chi_i)$. We denote
by $x_i$ a generator of $M(h_i,\chi_i)$,
whence $V=(x_1,\ldots,x_{\te})$ and $c$ is given
by the matrix $(b_{ij})$, with $b_{ij}=\chi_j(h_i)$.

\subsection{Lifting Procedure}\label{al:prl}
Let $A$ be a finite dimensional pointed Hopf algebra. Let
$A_0\subset A_1\subset\ldots$ be its coradical filtration, and
consider
$$H=\oplus_i\agg Hi=\gr A=\oplus_i(A_i/A_{i-1})\quad (A_{-1}=0).$$
It is a coradically graded Hopf algebra, i.e. $H_n=\oplus_{i=0}^n\agg Hi$.
Take $R=\coind H{\agg H0}$. This is a coradically graded Hopf algebra in
$\ydskg$, where $\Gm=G(A)$ is the group of group-likes of $A$
(and hence $\agg H0=\k\Gm$). Notice that the order of $\Gm$ divides
the dimension of $A$, by virtue of the Nichols--Zoeller theorem.
Furthermore, $\agg R0=\k 1$, from where $R$ satisfies \eqref{toba1},
\eqref{toba2} and \eqref{toba3}. If moreover $R$ satisfies \eqref{toba4}
then it is a Nichols algebra, and hence
$R=\toba{\agg R1}$. Since $R$ is an invariant of $A$,
we shall consider its ``infinitesimal part" defined by
\begin{equation}\label{eq:ddv}
V=V(A)=\agg R1,
\end{equation}
and call $\dim V$ {\it the rank of $A$}. It is clear that the
rank of $A$ is equal to $\dim A_1/\dim A_0-1$.

The lifting procedure
goes the other way: let $n\in\NN$ be fixed. We want to classify
the $n$-dimensional pointed Hopf algebras over $\k$. We can proceed
as follows: for each group $\Gm$ such that $|\Gm|$ divides $n$,
\begin{enumerate}
\item \label{pa:lp1}
	Determine the modules $V\in\ydskg$ such that
	$\dim\toba V\le \dfrac n{|\Gm|}$.
\item \label{pa:lp2}
	Classify the isomorphism classes of the bosonizations
	$\toba V\#\k\Gm$, where $V$ runs over the modules in step
	\ref{pa:lp1}. Note that two bosonizations $\toba V\#\k\Gm$ and
	$\toba V'\#\k\Gm$ may be isomorphic even though $\toba V$ and
	$\toba V'$ (thus, $V$ and $V'$) are not isomorphic (see \ref{lm:iso}).
\item \label{pa:lp3}
	For each isomorphism class $\toba V\#\k\Gm$ in step \ref{pa:lp2},
	classify all the liftings, i.e. all the pointed Hopf algebras
	(up to isomorphism) such that
	$\gr A=\toba V\#\k\Gm$.
\item \label{pa:lp4}
	Prove, if possible, that these are all the $n$-dimensional pointed
	Hopf algebras. Equivalently, show that a graded braided Hopf algebra
	$R$ in $\ydskg$ such that $\dim R=n/|\Gm|$ and satisfying
	\eqref{toba1}--\eqref{toba3}, also satisfies \eqref{toba4}
	(for this purpose it is useful to classify those $V$ in
	step \ref{pa:lp1} such that $\dim\toba V<n/|\Gm|$).
\end{enumerate}

\subsection{Contents}
The article is organized as follows: in section \ref{sprel} we compute
the Nichols algebras that appear applying the procedure above for
$32$-dimensional pointed Hopf algebras.

From \ref{sz2} to \ref{sgna} we follow the steps \ref{pa:lp1},
\ref{pa:lp2} and \ref{pa:lp3} of the procedure for each
of the groups of order $2^m$, $m\le 4$. In section \ref{srtoba}
we prove the statement in step \ref{pa:lp4} of the procedure.

\medskip
This work is part of the doctoral thesis of the author, who thanks
N. Andruskiewitsch for his excellent guidance.
I thank also H.-J. Schneider and M. Lamas for many valuable conversations.

Our main references on Hopf algebras are \cite{s} and \cite{mon}.
We refer to the survey article \cite{ag} for definitions and properties of
braided Hopf algebras, Nichols algebras (called there TOBAs) and for
further references on the subject.

\renewcommand{\theequation}{\arabic{section}.\arabic{equation}}
\section{The Nichols algebras we shall need} \label{sprel}

For future reference, we give two families of examples of Nichols
algebras. The first one corresponds to Taft algebras, and the second
one to the more general case of Quantum Linear Spaces (QLS for brevity).

\begin{exmp}\label{ex:taft}
Let $V$ be $1$-dimensional and let $x$ be a generator of $V$.
Suppose $c(x\otimes x)=qx\otimes x$, where $q\in\unid\k$.
If $q$ is a root of unity of order $n>1$, or if $q=1$ and $\car\k=n>0$,
then $\toba V=\k[x]/(x^n)$ as algebra (and hence $\dim\toba V=n$).
If $q$ is not a root of unity (which can happen for instance for
$V\in\yds{\k\ZZ}$), or if $q=1$ and $\car\k=0$, then $\toba V=\k[x]$ and
hence it is infinite dimensional. The comultiplication in both cases is
given by $\Dt(x^d)=\sum_{i=0}^d\ctr diq x^i\otimes x^{d-i}$.
\end{exmp}

\begin{rem}\label{rm:qn1}
The above example shows the following useful result: let $H$ be a
finite dimensional Hopf algebra and $R\in\yd$ a braided Hopf algebra
of dimension $N$. Suppose $\car\k$
does not divide $N$. If $x\in\prim R$ is such that $\k x$ is a
sub-YD-module of $R$ and $c(x\otimes x)=qx\otimes x$,
where $q\in\unid\k$, then $q\neq 1$.
\end{rem}
\begin{pf}
Suppose $q=1$. If $\car\k=0$ then $R$ would be
infinite dimensional, if $\car\k=p>0$ then $x$ would generate an
algebra of dimension a power of $p$, but then, taking the bosonization
of both algebras, $p$ would divide $N$.
\end{pf}

\begin{defn}\label{df:ddn}
Let $q\in\unid\k$. We define $N(q)$ by
$$N(q)=\begin{cases}
\mbox{order of $q$} & \mbox{if $q\neq 1$ and is a root of unity}, \\
\infty & \mbox{if $q$ is not a root of unity}, \\
\infty & \mbox{if $q=1$ and $\car\k=0$}, \\
\car\k & \mbox{if $q=1$ and $\car\k>0$}.
\end{cases}$$

\end{defn}

\begin{defn}
Let $V\in\yd$ be of dimension $\te$. We shall say that $V$ (or $c$)
{\it has a matrix} if there exists a basis $\{x_1,\ldots,x_\te\}$ of $V$
and a matrix
$(b_{ij})\in M_{\te}(\unid\k)$ such that the braiding $c$ is given by
$$c(x_i\otimes x_j)=b_{ij}x_j\otimes x_i.$$
When we want to specify the basis, we say that $V$ has a matrix
$(b_{ij})$ {\it in the basis} $\{x_1,\ldots,x_\te\}$.
If $H=\k\Gm$ and $\Gm$ is abelian, then every finite dimensional
module has a matrix (see \eqref{eq:bca}).

Let $V\in\ydskg$ be finite dimensional.
We shall say that $V$ (or $c$) {\it comes from the abelian case} if
there exists an invariant abelian subgroup $\Gm'\subset\Gm$ such that
the image of $\dt:V\to\k\Gm\otimes V$ is included in $\k\Gm'\otimes V$.
In this case, $V$ can be considered as a module in $\yds{\k\Gm'}$
(thus, it has a matrix) and the bosonization $\toba V\#\k\Gm$ can be
reconstructed as an extension of $\Gm/\Gm'$ by $\toba V\#\k\Gm'$.
We note that the condition on $\Gm'$ to be invariant can be dropped:
it is not difficult to see that if there exists $\Gm'$ abelian such that
the image of $\dt$ lies in $\k\Gm'\otimes V$, then there exists
$\Gm''\subset\Gm'$ invariant (and clearly abelian) such that the image
of $\dt$ lies in $\k\Gm''\otimes V$.
\end{defn}

\begin{exmp}\label{ex:qls}
Let $V\in\yd$ be of dimension $\te$, and suppose $V$ has a matrix
$(b_{ij})$ in the basis $\{x_1,\ldots,x_\te\}$.
Let $N_i=N(b_{ii})$, see \ref{df:ddn}.
Then it can be proved (see \cite[Prop. 3.5]{as1}) that
$$\dim\toba V\ge\prod_{i=1}^nN_i.$$
Furthermore, if $N_i<\infty$ for all $i$, then
the equality holds if and only if
\begin{equation}\label{eq:qls}
b_{ij}b_{ji}=1\mbox{ for all } i\neq j.
\end{equation}
In this case, a PBW basis of $\toba V$ is given by
$$\{x_1^{n_1}\cdots x_{\te}^{n_{\te}}\ |\ 0\le n_i<N_i\},$$
and the multiplication is given by the following relations:
$$x_i^{N_i}=0,\quad x_jx_i=b_{ji}x_ix_j\mbox{ for }i<j.$$
Hence, as an algebra, $\toba V$ is a quantum linear space.
Let now $H=\k\Gm$. In order to agree with the notation of \cite{as1},
we shall say that $\toba V$ is a {\it Quantum Linear Space} (or QLS)
if $V$ satisfies \eqref{eq:qls} and the lines $\k x_i$ are submodules of
$V$ in $\ydskg$. If $V$ can be considered as a module in $\yds{\k\Gm'}$ for
certain subgroup $\Gm'\subset\Gm$ (that is, if the image of the coaction
$\dt$ lies in $\Gm'\otimes V$) and it is a QLS in $\yds{\k\Gm'}$ then
we say that it is a QLS {\it over} $\Gm'$.
\end{exmp}

\bigskip
We shall use an important result about Nichols algebras due to
Nichols:
\begin{prop}\label{pr:dualpuan}
Let $R=\oplus_{i=0}^N\agg Ri$ be a Nichols algebra (in some braided
rigid abelian category). Suppose that $\agg RN\neq 0$. Then
$$\dim\agg Ri=\dim\agg R{N-i}\quad\forall i=0,\ldots,N.$$
\end{prop}
\begin{pf}
See for instance \cite[Pr. 3.2.2]{ag}.
\end{pf}

\begin{lem}\label{lm:nht}
Let $A$ be a $2^n$-dimensional pointed Hopf algebra with $A_0=\k\Gm$.
Suppose $|\Gm|\ge 2^{n-2}$ and let $V=V(A)\in\ydskg$ be the module
defined by \eqref{eq:ddv}. Then $\dim V\le 2$.
\end{lem}
\begin{pf}
Since $2^n=|\Gm|\times\dim\toba V$ we have $\dim\toba V\le 4$;
but if $\dim V\ge 3$, then $\dim\toba V\ge 5$ by \ref{pr:dualpuan}.
\end{pf}

In \cite{n} the author introduces derivations for Nichols algebras,
which turns out to be a useful tool for computing bases (and hence their
dimensions).
\begin{defn}\label{df:der}
Let $V\in\yd$. For $i+j=n$, we denote by
$$\Dt_{i,j}:\tobag nV\to\tobag iV\otimes\tobag jV$$
the $(i,j)$-component of the comultiplication of $\toba V$.
It is proved in \cite{schb} (or see \cite[Def. 3.2.10]{ag})
that $\Dt_{i,j}$ is injective $\forall i,j$. Let
$\{x_1,\ldots,x_\te\}$ be a basis of $V$ and let
$\{x^*_1,\ldots,x^*_\te\}$ be its dual basis. We denote by $\ptl_{x_i}$
the differential operator on $\toba V$ given by
$$\ptl_{x_i}(z)=(\id\otimes x^*_i)\Dt_{n-1,1}(z),
	\quad\mbox{if } z\in \tobag nV,\ n>0,
	\quad\mbox{and }\ptl_{x_i}(1)=0.$$
By the injectivity of $\Dt_{i,j}$ it is clear that for $z\in\tobag nV$
($n>0$) we have $z=0$ if and only if
$\ptl_{x_i}(z)=0$ for all $i=1,\ldots,\te$.
Suppose now that $V\in\ydskg$ and some $x^*_i$ vanishes on the homogeneous
components of degree $\neq g_i$ for some $g_i\in\Gm$ (that is, $x^*_i(y)=0$
if $\dt(y)=h\otimes y$ and $h\neq g_i$).
Then it is easy to see that $\ptl_{x_i}$ verifies the Leibniz rule
$$\ptl_{x_i}(z_1z_2)=\ptl_{x_i}(z_1)(g_i\ganda z_2)+z_1\ptl_{x_i}(z_2).$$
Since usually we shall take bases $\{x_1,\ldots,x_\te\}$ of homogeneous
elements (i.e.
$\dt(x_i)=g_i\otimes x_i$ for some $g_i\in\Gm$), the previous
condition will hold for all $x^*_i$.
\end{defn}

\begin{defn}
Let $A$ be a Hopf algebra. We shall denote by $\Ad$ the
adjoint action
$$\Ad_x(y)=\com x1y\Ss(\com x2).$$
Let $B$ be a braided Hopf algebra with multiplication $\mu$,
let $c\in\End(B\otimes B)$ be the braiding. We shall denote by $\Ad$
the braided adjoint action, given by
$$\Ad_x(y)=\mu(\mu\otimes\Ss)(\id\otimes c)(\Dt\otimes\id)(x\otimes y).$$
Note that if $x\in\prim B$ then the adjoint is nothing but
$\Ad_x(y)=\mu(\id-c)(x\otimes y)$.
Since in the cases we are interested in the endomorphism $c$ is given
by the context, we do not make a reference to $c$ in the notation.
\end{defn}

We need the following relation between both adjoint actions:
let $H$ be a graded Hopf algebra, $\pi:H\to\agg H0$ the canonical projection
and $B\in\yds{\agg H0}$ the braided Hopf algebra of coinvariants (then
$H\simeq B\#\agg H0$). Let $x\in\agg B1$ (it is a primitive element).
Then it is straightforward to see that
$$\Ad_{(b\# 1)}(c\# 1)=(\Ad_b(c))\# 1\quad\forall c\in B,$$
where the left $\Ad$ is the classical one, while the right $\Ad$ is the
braided one.

\bigskip
We shall compute now the Nichols algebras appearing in next sections.
\begin{prop}\label{lm:lsn}
Let $V\in\yd$ with basis $\{x,y\}$. Suppose that $c$ has a matrix
$(b_{ij})$ in this basis, and suppose that $b_{12}b_{21}\neq 1$,
$b_{11}\neq 1$, $b_{22}\neq 1$. Let $z=\Ad_x(y)=xy-b_{12}yx$.
Let $N_i=N(b_{ii})$ (see \ref{df:ddn}) and suppose that $N_i<\infty$
(i.e., that $b_{ii}$ is a root of unity, or $b_{ii}=1$ and $\car\k>0$).
Let $q=b_{11}b_{12}b_{21}b_{22}$.
\begin{enumerate}
\item 
There exists $n_{12}$ (respectively $n_{21}$) such that
$(\Ad_x)^{n_{12}}(y)=0$ (respectively
$(\Ad_y)^{n_{21}}(x)=0$). Let $n_{12}$ (resp. $n_{21}$) be the
least positive integer with such property. Then $n_{12}\le N_1$
and $n_{21}\le N_2$. Suppose that there exists $m_{12}$
(resp. $m_{21}$) such that $b_{12}b_{21}b_{11}^{m_{12}}=1$
(resp. $b_{12}b_{21}b_{22}^{m_{21}}=1$) and let $m_{12}$ (resp. $m_{21}$)
be the least non negative integer with such property.
Then $n_{12}=m_{12}+1$ (resp. $n_{21}=m_{21}+1$).
\item If $q=1$ and $\car\k=0$ then $z^n\neq 0\ \forall n$ and then
$\toba V$ is infinite dimensional. If $q=1$ and $\car\k\neq 0$ then
$\dim\toba V\ge N_1N_2\car\k$.
\item Suppose $b_{11}=b_{22}=-1$ and let $N=N(q)$.
Then $\toba V$ has dimension $4N$, with a PBW
basis $\{x^ay^bz^c\ |\ 0\le a,b<2,\ 0\le c<N\}$.
\end{enumerate}
\end{prop}

\begin{pf}
\begin{enumerate}
\item We compute inductively the derivatives of $(\Ad_x)^n(y)$
(see \ref{df:der}):
$$\ptl_x((\Ad_x)^1(y))=\ptl_x(xy-b_{12}yx)=0.$$
Suppose that $\ptl_x((\Ad_x)^n(y))=0$, then
\begin{align*}
\ptl_x((\Ad_x)^{n+1}(y))&=
	\ptl_x\left(x(\Ad_x)^n(y)-b_{11}^nb_{12}(\Ad_x)^n(y)x\right) \\
&=b_{11}^nb_{12}(\Ad_x)^n(y)-b_{11}^nb_{12}(\Ad_x)^n(y)=0, \\
\ptl_y((\Ad_x)^1(y))&=\ptl_y(xy-b_{12}yx)=(1-b_{12}b_{21})x.
\end{align*}
Let $\sgm_n=(1-b_{12}b_{21}b_{11}^n)$.
Suppose that $\ptl_y((\Ad_x)^n(y))=\sgm_0\sgm_1\cdots\sgm_{n-1}x^n$. Then
\begin{align*}
\ptl_y((\Ad_x)^{n+1}(y))&=\ptl_y(x(\Ad_x)^n(y)-b_{11}^nb_{12}(\Ad_x)^n(y)x) \\
&=\sgm_0\cdots\sgm_{n-1}x^{n+1}-b_{11}^nb_{12}\sgm_0\cdots\sgm_{n-1}b_{21}x^{n+1} \\
&=\sgm_0\cdots\sgm_nx^{n+1}.
\end{align*}
Thus, since $x^{N_1}=0$, we have
$$\ptl_x((\Ad_x)^{N_1}(y))=\ptl_y((\Ad_x)^{N_1}(y))=0,$$
and hence $(\Ad_x)^{N_1}(y)=0$, which proves that $n_{12}$ exists
and $n_{12}\le N_1$.
Furthermore, if $b_{12}b_{21}b_{11}^{m_{12}}=1$, we see that
$$\ptl_x((\Ad_x)^{m_{12}+1}(y))=\ptl_y((\Ad_x)^{m_{12}+1}(y))=0,$$
whence $n_{12}\le m_{12}+1$ (notice that if $m_{12}$ exists, then
$m_{12}<N_1$). Moreover, the elements
$x^i$ do not vanish for $i<N_1$, whence if $m_{12}$ exists then
$\ptl_y((\Ad_x)^{m_{12}}(y))\neq 0$, from where $n_{12}>m_{12}$ and thus
$n_{12}$ is exactly $m_{12}+1$.

\item \label{abo} We compute the derivatives of $z^n$. We do it for a
general $q=b_{11}b_{12}b_{21}b_{22}$ (keep the notation
$\sgm_n=(1-b_{12}b_{21}b_{11}^n)$).
\begin{align*}
\ptl_x(z)&=\ptl_x(\Ad_xy)=0, \\
\ptl_y(z)&=\ptl_y(\Ad_xy)=\sgm_0x, \\
\ptl_x\ptl_y(z)&=\ptl_x(\sgm_0x)=\sgm_0,
\end{align*}
\begin{align*}
\ptl_y(z^n)&=\sum_{i=0}^{n-1}(b_{21}b_{22})^i\sgm_0z^{n-1-i}xz^i, \\
\ptl_x\ptl_y(z^n)
	&=\sum_{i=0}^{n-1}(b_{21}b_{22})^i(b_{11}b_{12})^i\sgm_0z^{n-1} \\
	&=(n)_q\sgm_0z^{n-1}.
\end{align*}
Let now $q=1$. If $\car\k=0$, we see inductively that $z^n\neq 0$
$\forall n$. Furthermore, since $z^n\in\tobag{2n}V$, the set
$\{z^n\ |\ n\ge 0\}$ is linearly independent, and thus $\toba V$
is infinite dimensional.
If $\car\k>0$, we see inductively that $z^n\neq 0$ $\forall n<\car\k$.
Furthermore, with the same computations as in \ref{bel} below, we see
that the set
$$\{x^ay^bz^c\ |\ 0\le a<N_1,\ 0\le b<N_2,\ 0\le c<\car\k\}$$
is linearly independent, whence $\dim\toba V\ge N_1N_2\car\k$.

\item \label{bel} By the same computations than in \ref{abo}, we have that
the set $\{z^i\ |\ 0\le i<N\}$ is linearly independent. We first prove that
$z^N=0$:
\begin{align*}
\ptl_x(z^N)&=0 \quad\mbox{(as in \ref{abo})} \\
\ptl_y(z^N)&=\sum_{i=0}^{N-1}(b_{21}b_{22})^i\sgm_0z^{N-1-i}xz^i,
\end{align*}
but this element vanishes in $\toba V$, since
\begin{align*}
\ptl_x(\ptl_yz^N)&=(N)_q\sgm_0z^{N-1}=0\quad\mbox{(as in \ref{abo})} \\
\ptl_y(\ptl_yz^N)&=\sum_{i=0}^{N-1}(b_{21}b_{22})^i\sgm_0
	[\sum_{j=0}^{N-2-i}(b_{21}b_{22})^{j+i}b_{21}z^{N-2-i-j}xz^jxz^i \\
&\hspace{2cm} +\sum_{j=0}^{i-1}(b_{21}b_{22})^jz^{N-1-i}xz^{i-j-1}xz^j] \\
&=\sum_{k+l\le N-2}(b_{21}b_{22})^{k+2l}\sgm_0b_{21}(1+b_{22})
		z^{N-2-k-l}xz^kxz^l \\
&=0\mbox{ (since $b_{22}=-1$)}.
\end{align*}
Now, since $b_{11}=b_{22}=-1$, we have $x^2=y^2=0$. Furthermore,
\begin{align*}
yx&=b_{12}^{-1}xy-b_{12}^{-1}z, \\
zx&=(xy-b_{12}yx)x=xyx=-b_{12}^{-1}x(xy-b_{12}yx)=-b_{12}^{-1}xz, \\
zy&=(xy-b_{12}yx)y=-b_{12}yxy=-b_{12}y(xy-b_{12}yx)=-b_{12}yz.
\end{align*}
Thus the set $\{x^ay^bz^c\ |\ 0\le a,b<2,\ 0\le c<N\}$ generates
$\toba V$ as a vector space. Let $\al_{b,c}\in\k$ such that
$$\sum\begin{Sb}0\le b<2 \\ 0\le c<N\end{Sb} \al_{b,c}y^bz^c=0.$$
We compute
\begin{align*}
\ptl_y(yz^c)&=(b_{21}b_{22})^cz^c+y\ptl_yz^c, \\
\ptl_x\ptl_y(y)&=0, \\
\ptl_x\ptl_y(yz^c)&=(c)_q\sgm_0yz^{c-1}\quad(c>0).
\end{align*}
Applying $(\ptl_x\ptl_y)^{N-1}$, we see that $\al_{0,N-1}=\al_{1,N-1}=0$.
Inductively, we see that $\al_{b,c}=0\ \forall b,c$.
In the same vein, let $\al_{a,b,c}\in\k$ such that
$$\sum\begin{Sb}0\le a,b<2 \\ 0\le c<N\end{Sb} \al_{a,b,c}x^ay^bz^c=0.$$
We compute
$$\ptl_x(xy^bz^c)=b_{11}^cb_{22}^{b+c}y^bz^c,$$
from where we see that $\al_{1,b,c}=0\ \forall b,c$,
and thus $\al_{a,b,c}=0\ \forall a,b,c$, which proves the statement.
\end{enumerate}
\ 
\end{pf}

\begin{rem}
It can be proved in general that if $V$ has a matrix $(b_{ij})$
in a basis $\{x_1,\ldots,x_\te\}$, then the operator
$(\ptl_{x_i})^{N(b_{ii})}$ vanishes.
Furthermore, the duality of \cite[Prop. 3.2.30]{ag} translates
products into compositions.
\end{rem}

\begin{defn}
Let $V$ be a $\te$-dimensional YD module, and suppose
$V$ has a matrix $(b_{ij})$. According to \cite{as2}, we say that
$V$ {\it is of Cartan type} (or CT) if there exists a matrix
$(a_{ij})\in M_\te(\ZZ)$ such that
$$b_{ij}b_{ji}=b_{ii}^{a_{ij}}\quad\forall i,j.$$
If we choose $a_{ii}=2$ and $-\mbox{order of }b_{ii}<a_{ij}\le 0$ for
$i\neq j$, then $(a_{ij})$ is a generalized Cartan matrix. We transfer
to $V$ the terminology on $(a_{ij})$.
\end{defn}

\begin{cor} \label{pr:cto}
\hspace*{1in}
\begin{enumerate}
\item
Let $V\in\yd$, and suppose that $V$ has one of the following matrices
in the basis $\{x,y\}$ (recall from the beginning of section \ref{ss:bdac}
that $\rmu,\xi$ denote fixed roots of unity of orders $4$ and $8$ respectively).
$$b^1_q=\mdpd {-1}{\fsm q^{-1}}{q}{-1^{\fmu}}\ (q\in\unid\k),
	\qquad b^2=\mdpd {-1}{-1}{1}{-1},$$
$$b^3_{\pm+}=\mdpd {-1}{1}{\pm \rmu}{-1},\qquad
	b^3_{\pm-}=\mdpd {-1}{-1}{\pm \rmu}{-1},$$
$$b^4_{\pm}=\mdpd {\pm \rmu}{\pm \rmu}{-1}{-1},\qquad
	b^5_{\pm}=\mdpd {\pm \rmu}{\mp \rmu}{-1}{-1},$$
$$b^6_k=\mdpd {-1}{1}{\xi^k}{-1}\quad (k=1,3,5,7).$$
The information about $\toba V$ that we need is given in the following
table. We put $z_1=\Ad_x(y)$ or $z_1=\Ad_y(x)$ depending on the case.
$$\begin{array}{|l|c|c|c|r|r|r|} \hline
\multicolumn{7}{|c|}{\mbox{TABLE 1:
	Matrices and algebras}} \\ \hline
\mbox{b} & \dim\toba V & \mbox{CT?} & z_1
	& |x| & |y| & |z_1|  \\ \hline
b^1_q           & 4     & \mbox{Yes} & 0 & 2 & 2 &  \\ \hline
b^2             & 8     & \mbox{Yes} & \Ad_x(y)=xy+yx
	& 2 & 2 & 2	\\ \hline
b^3_{\pm+}      & 16	& \mbox{No} & \Ad_x(y)=xy+yx
	& 2 & 2 & 4	\\ \hline
b^3_{\pm-}      & 16    & \mbox{No} & \Ad_x(y)=xy+yx
	& 2 & 2 & 4	\\ \hline
b^4_{\pm}       & 16	& \mbox{No} & \Ad_y(x)=xy+yx
	& 4 & 2 & 2	\\ \hline
b^5_{\pm}       & \ge 8N(1)\ge 24& \mbox{No} & \Ad_y(x)=xy+yx
	& 4 & 2 & N(1)	\\ \hline
b^6_k           & 32    & \mbox{No} & \Ad_y(x)=xy-yx
	& 2 & 2 & 8	\\ \hline
\end{array}$$
We denote by $|w|$ the nilpotency order of $w$, i.e.
$w^{|w|}=0,\ w^{|w|-1}\neq 0$.

\item
Suppose $V$ is $3$-dimensional and has a matrix $b^1_{q_1,q_2,q_3}$
in the basis $\{x,y,z\}$, where
$$b^1_{q_1,q_2,q_3}=\left(\begin{array}{rrr}
	-1^{\fmu} &q_1^{-1}&q_3^{-1} \\
	q_1^{\fmu} &-1^{\fmu} &q_2^{-1} \\
	q_3^{\fmu} &q_2^{\fmu} &-1^{\fmu} 
\end{array}\right)\quad (q_1,q_2,q_3\in\unid\k).$$
Hence $\toba V$ is an $8$-dimensional QLS, with relations
$$x^2=y^2=z^2=0,\qquad yx=q_1xy,\ zx=q_3xz,\ zy=q_2yz.$$
\end{enumerate}
\end{cor}

\begin{pf}
The cases $b^1_q$, $b^1_{q_1,q_2,q_3}$ are particular cases of
example \ref{ex:qls}.

The cases $b^2$, $b^3_{\pm+}$, $b^3_{\pm-}$, $b^5_{\pm}$ and $b^6_k$
are particular cases of \ref{lm:lsn}.

For $b^4_{\pm}$, it is straightforward to see that
$$zy=yz,\qquad zx=-\rmu xz,$$
and, using the same techniques than in \ref{lm:lsn}, that there is
a PBW basis given by \\
$\{x^ay^bz^c\ |\ 0\le a<4,\ 0\le b,c<2\}$.
\end{pf}

\begin{rem}
Let $V$ be an object in $\cate$, a braided abelian category.
Suppose $V$ has a matrix
$$b=\mdpd {b_{11}}{b_{12}}{b_{21}}{b_{22}}$$
in the basis $\{x,y\}$. Then it is clear that $V$ has the matrix
$$[b]^{\rho}:=\mdpd {b_{22}}{b_{21}}{b_{12}}{b_{11}}$$
in the basis $\{y,x\}$.

Consider now $\cateb$ the category with same objects and morphisms
as $\cate$ but with the inverse braiding $\bar c=c^{-1}$. Consider
$\bar V=V$ in $\cateb$. Then $\bar V$ has a matrix
$$[b]^{\tau}:=\mdpd {b_{11}^{-1}}{b_{21}^{-1}}{b_{12}^{-1}}{b_{22}^{-1}}$$
in the basis $\{x,y\}$. On the other hand, it can be shown that 
$\dim\toba{\bar V}=\dim\toba V$.
Notice, for instance, that $[b^3_{\pm +}]^{\rho\tau}=b^3_{\mp +}$.
\end{rem}

Let $M\in\ydskg$ and $f\in\Aut(\Gm)$. We write $(M)^f$ for the YD-module
with the same underlying space as $M$ but with the structure given by
$$h\ganda^fm=f^{-1}(h)\ganda m,\quad \dt^f(m)=(f\otimes\id)\dt(m).$$
Note that $M(g,\rho)^f=M(f(g),\rho\circ f^{-1})$.
In order to classify non isomorphic algebras, the following criterium will
be applied:

\begin{lem}\label{lm:iso}
Let $V=\oplus_{i=1}^nM(g_i,\rho_i)$ and
$V'=\oplus_{i=1}^{n'}M(g'_i,\rho'_i)$,
where $\rho_i$ and $\rho'_i$ are irreducible.
Then $\toba V\#\k\Gm\simeq\toba {V'}\#\k\Gm$ if and only if
$n=n'$ and there exist $f\in\Aut(\Gm)$, $\sgm\in\SS_n$ such that
$M(g'_i,\rho'_i)\simeq(M(g_{\sgm(i)},\rho_{\sgm(i)}))^f$.
\end{lem}

\begin{pf}
Mimic \cite[Prop. 6.3]{as2}, where it is proved for $\Gm$ abelian.
\end{pf}

\bigskip
For the last step of the lifting procedure we shall use the following
Proposition, which gives an invariant way to lift the elements of the
bosonized algebras.

\begin{prop}
Let $A$ be a pointed Hopf algebra, $H=\gr A$ be its 
associated graded Hopf algebra and $\Gm$ the group of group-likes of $A$.
Let 
$V=V(A)\in\ydskg$ be the YD-module associated to $A$ as in \eqref{eq:ddv},
that is, $V=\coind {(\agg H1)}{\agg H0}$. We have
$\agg H1=V\#\k\Gm$. Take on $\agg H1$ the $\Gm$ action given by
$$g\cdot(v\# h)=(g\ganda v)\# ghg^{-1}$$
Consider on $A$ the adjoint action of $\Gm$. It restricts to the
coradical filtration, and then $A_1$ becomes a $\k\Gm$-module.
Let $p_0:A_1\to\agg H1=A_1/A_0$ be the canonical projection. Then
\begin{enumerate}
\item $p_0$ is a morphism of $\k\Gm$-modules. \label{p01}
\item Let $\lmb\in\k^{\Gm}$, $\lmb(g)=\dt_{g,1}$. Then
	$p=(\id\#\lmb)\circ p_0:A_1\to V$ is a morphism
	of $\k\Gm$-modules. \label{p02}
\item There exists a section $s:V\to A_1$ of $\k\Gm$-modules to $p$
	such that if 
	$\dt(x)=g\otimes x$ then $s(x)\in\primg Ag$.\label{p03}
\end{enumerate}
\end{prop}

\begin{pf}
The items \ref{p01} and \ref{p02} are straightforward looking at
the definition of the action $\ganda$. We prove  \ref{p03}.
For $g,h\in\Gm$, let $\primghp gh\subset A_1$ be subspaces such that
$\primgh Agh=\k(g-h)\oplus\primghp gh$. Then, by the Taft--Wilson theorem,
$$A_1=\k\Gm\oplus(\bigoplus_{g,h\in\Gm}\primghp gh).$$
We have the exact sequence of vector spaces
$$0\to\k\Gm=A_0\to A_1\to V\#\k\Gm\to 0.$$

Consider now the inclusion
$$i:\oplus_{g\in\Gm}\primgp g\inclu A_1.$$
It is easy to verify that the composition $p\circ i$ is an isomorphism of
vector spaces (a word of warning is needed here: $\oplus_{g\in\Gm}\primgp g$
need not to be a $\k\Gm$-submodule of $A_1$). Let
$s_0:V\to \oplus_{g\in\Gm}\primgp g$, $s_0=(p\circ i)^{-1}$, and consider
as usual
$$s:V\to A_1,\quad s=\frac 1{|\Gm|}\sum_{h\in\Gm}(h\cdot s_0),
	\mbox{ i.e. }
	s(x)=\frac 1{|\Gm|}\sum_{h\in\Gm}h\cdot(s_0(h^{-1}\ganda x)).$$
It is clear that $s$ is a $\k\Gm$ section to $p$.

Let $V^g=\{x\in V\ |\ \dt(x)=g\otimes x\}$; then $V=\oplus_{g\in\Gm}V^g$.
Moreover, let $a\in\primgp g$. Then $\Dt a=g\otimes a+a\otimes 1$, and
if $b=a+\k\Gm\in A_1/A_0$ is the image of $a$ under $A_1\to A_1/A_0$, we
have $\Dt b=g\otimes b+b\otimes 1$. Then $b\in\coind{(V\#\k\Gm)}{\agg H0}$,
whence $b=v\# 1$ for some $v\in V$ such that
$\Dt(v\# 1)=g\otimes (v\# 1)+(v\# 1)\otimes 1$.
This implies that $\dt(v)=g\otimes v$, and hence
$v=(\id\#\lmb)(v\# 1)\in V^g$. Thus $p(\primgp g)\subseteq V^g$
for all $g\in\Gm$. Now, since
$p\circ i:\oplus_{g\in\Gm}\primgp g\to \oplus_{g\in\Gm}V^g$ is an
isomorphism, we must have that
$p\circ i|_{\primgp g}:\primgp g\to V^g$ is an isomorphism,
whence $s_0(V^g)=\primgp g$ for all $g\in\Gm$.

Let now $x\in V^g$. We claim that $s(x)\in\primg Ag$. This is true for,
by the definition of $s$, we have
\begin{align*}
s(x)&=\frac 1{|\Gm|}\sum_{h\in\Gm}h\cdot s_0(h^{-1}\ganda x)
	\in\sum_{h\in\Gm}h\cdot s_0(V^{h^{-1}gh}) 
	=\sum_{h\in\Gm}h\cdot\primgp{h^{-1}gh} \\
& \subseteq \sum_{h\in\Gm}h\cdot\primg A{h^{-1}gh} 
	=\sum_{h\in\Gm}\primg Ag=\primg Ag.
\end{align*}
\end{pf}

We will use the above proposition as follows:
\begin{tactic}\label{ta:lev}
If $V$ has a basis $\{x_1,\ldots,x_\te\}$ such that $x_i\in V^{g_i}$,
the following proposition enables us to
consider $a_1,\ldots,a_\te$ such that $a_i\in\primg A{g_i}$, the
projection $A_1\to A_1/A_0$ sends $a_i$ to $(x_i\# 1)\in V\#\k\Gm$,
and if $g\in\Gm$ acts
 on $V$ by some matrix $T$ in the basis $\{x_i\}$,
then $g$ acts on the linear span $(a_1,\ldots,a_\te)$ by the same matrix
(in the basis $\{a_i\}$).
Furthermore, since $V$ generates $\toba V$, the $a_i$'s together with
$\Gm$ generate $A$. In general, suppose
$\sum_i(\al_ix_{i,1}\cdots x_{i,n})=0\in \tobag nV$.
Then $b=\sum_i(\al_ia_{i,1}\cdots a_{i,n})\in A$ lies
in the kernel of the projection, i.e. $b\in A_{n-1}$. Usually, we
shall have $b\in\primg Ag$ for some $g\in\Gm$, from where $b$ shall lie
in $A_1$.
\end{tactic}

\smallskip
The following is a very easy but useful remark:
\begin{rem}\label{rm:ebu}
Let $A$ be a finite dimensional Hopf algebra such that $\car\k$
does not divide $\dim A$.
Let $g\in A$ be a group-like and $a\in\primg Ag$ such that $ga=\zeta ag$.
Let $n$ be the order of $\zeta$. Then, by the same computations as in
\ref{ex:taft}, $a^n\in\primg A{g^n}$. If $g^n=1$ we must have $a^n=0$,
for, if not, $a$ would generate a sub Hopf algebra of dimension a power
of $\car\k$ (if $\car\k>0$) or $\infty$ (if $\car\k=0$), leading to a
contradiction.

If we are lifting a QLS, we have $a_i\in\primg A{g_i}$ and
$g_ia_jg_i^{-1}=b_{ij}a_j$.
If $i\neq j$, $b_{ij}b_{ji}=1$ and $g_ig_j=g_jg_i$. Then
\begin{align*}
\Dt(a_ia_j-b_{ij}a_ja_i) &= g_ig_j\otimes a_ia_j+g_ia_j\otimes a_i
	+a_ig_j\otimes a_j+a_ia_j\otimes 1 \\
	&\quad -b_{ij}g_jg_i\otimes a_ja_i
	-b_{ij}a_jg_i\otimes a_i-b_{ij}g_ja_i\otimes a_j
	-b_{ij}a_ja_i\otimes 1 \\
&=g_ig_j\otimes (a_ia_j-b_{ij}a_ja_i)+(a_ia_j-b_{ij}a_ja_i)\otimes 1,
\end{align*}
from where $a_ia_j-b_{ij}a_ja_i\in\primg A{g_ig_j}$.
As above, if $\car\k\not|\dim A$ and $g_ig_j=1$, then
$a_ia_j-b_{ij}a_ja_i=0$.
\end{rem}

\section{$\Gm=C_2$} \label{sz2}
Let $V=\oplus_iM(h_i,\chi_i)$. For $\toba V$ to have dimension $16$,
we must have 
$b_{ii}=-1\ \forall i$, where $(b_{ij})=(\chi_j(h_i))$
is the matrix of $V$. Thus, $V=\oplus_iM(g,\hat g)$, and hence
$(b_{ij})=\hat g(g)=-1\ \forall ij$, from where $\toba V$ is a QLS with
generators $x_1,\ldots,x_{\te}$ and relations $x_i^2=0$, $x_ix_j=-x_jx_i$.
Thus, $\dim\toba V=2^\te$, and $\te=4$. Let $a_1,\ldots,a_4$ be as
in \ref{ta:lev}. Then $ga_ig=-a_i$ and $a_i\in\primg Ag$. We have
by \ref{rm:ebu} that $a_i^2\in\prim A$, and the same is true
for $a_ia_j+a_ja_i$. Hence $a_i^2=a_ia_j+a_ja_i=0$ $\forall i,j$.
There exist thus only one Hopf
algebra of dimension $32$ with coradical $\k C_2$.
It is clear, using the same reasoning, that there exists a Hopf
algebra with coradical $\k C_2$ and dimension $n$ if and only if
$n=2^m$ for some $m\in\NN$. Furthermore, if $n=2^m$, there exists only
one such Hopf algebra. This is the content of
\cite[Thm. 4.2.1]{n}, but we included here a proof for completeness.

\renewcommand{\theequation}{\arabic{section}.\arabic{subsection}.\arabic{equation}}
\section{$\Gm$ with order $4$}
\subsection{$\Gm=C_2\times C_2$} \label{sz2z2}
Necessarily the rank is $2$ or $3$: for rank $1$ the root $b_{11}$ has
order $2$, whence $\dim(\toba V\#\k\Gm)$ would be $8$, and rank $\ge 4$
is impossible by \ref{ex:qls}.
The following is a list of the irreducible \ydm s giving Nichols
algebras of dimension $2$, with the
action of $\Aut(\Gm)$ as in the lemma \ref{lm:iso}.
In this case, $\Aut(\Gm)$ is isomorphic to $\SS_3$, generated by
$f_1$ and $f_2$, where, using the notation at the end of \ref{ss:bdac},
$$f_1=(g_1\mapsto g_1g_2,\ g_2\mapsto g_2),\quad
f_2=(g_1\mapsto g_1,\ g_2\mapsto g_1g_2).$$
Thus, we have
$$\begin{array}{|c|l|l|c|c|} \hline
\multicolumn{5}{|c|}{\mbox{TABLE 2:
	modules $/\Aut(\Gm)$, $\Gm=C_2\times C_2$}} \\ \hline
M(h,\chi) & h & \chi & (M(h,\chi))^{f_1} & (M(h,\chi))^{f_2} \\ \hline
Y_2^1 & g_1 & \hat g_1          & Y_2^5 & Y_2^2 \\ \hline
Y_2^2 & g_1 & \hat g_1\hat g_2  & Y_2^6 & Y_2^1 \\ \hline
Y_2^3 & g_2 & \hat g_2          & Y_2^4 & Y_2^6 \\ \hline
Y_2^4 & g_2 & \hat g_1\hat g_2  & Y_2^3 & Y_2^5 \\ \hline
Y_2^5 & g_1g_2 & \hat g_1       & Y_2^1 & Y_2^4 \\ \hline
Y_2^6 & g_1g_2 & \hat g_2       & Y_2^2 & Y_2^3 \\ \hline
\end{array}$$

Consider then $V=M_1\oplus M_2=M(h_1,\chi_1)\oplus M(h_2,\chi_2)$.
Since $\Aut(\Gm)$ acts transitively, we may suppose that $M_1=Y_2^1$.
For rank $2$ this gives $6$ possibilities, but the cases $Y_2^1\oplus Y_2^4$
and $Y_2^1\oplus Y_2^6$ are equivalent by means of $f_1\circ f_2$. Hence, the
possibilities are

$$\begin{array}{|c|l|l|l|l|l|c|} \hline
\multicolumn{7}{|c|}{\mbox{TABLE 3:
	Rank $2$, $\Gm=C_2\times C_2$}} \\ \hline
V     & h_1 & h_2 & \chi_1 & \chi_2 & (b_{ij}) & \dim\toba V \\ \hline
V_2^1 & g_1 & g_1 & \hat g_1 & \hat g_1           & b^1_{-1} & 4 \\ \hline
V_2^2 & g_1 & g_1 & \hat g_1 & \hat g_1\hat g_2   & b^1_{-1} & 4 \\ \hline
V_2^3 & g_1 & g_2 & \hat g_1 & \hat g_2           & b^1_1    & 4 \\ \hline
V_2^4 & g_1 & g_2 & \hat g_1 & \hat g_1\hat g_2   & b^2      & 8 \\ \hline
V_2^5 & g_1 & g_1g_2 & \hat g_1 & \hat g_1        & b^1_{-1} & 4 \\ \hline
\end{array}$$

\medskip
For rank $3$, $V=\oplus_{i=1}^3M_i=\oplus_{i=1}^3M(h_i,\chi_i)$.
Since the module $V_2^4$ must not appear as a submodule, we have
the following possibilities:
$$\begin{array}{|c|l|l|l|l|l|l|c|c|} \hline
\multicolumn{9}{|c|}{\mbox{TABLE 4:
	Rank $3$, $\Gm=C_2\times C_2$}} \\ \hline
V      & h_1 & h_2 & h_3 & \chi_1 & \chi_2 & \chi_3
		& (b_{ij}) & \dim\toba V \\ \hline
W_2^1 & g_1 & g_1 & g_1 & \hat g_1 & \hat g_1 & \hat g_1 
		& b^1_{-1,-1,-1} & 8 \\ \hline
W_2^2 & g_1 & g_1 & g_1 & \hat g_1 & \hat g_1 & \hat g_1\hat g_2 
		& b^1_{-1,-1,-1} & 8 \\ \hline
W_2^3 & g_1 & g_1 & g_1g_2 & \hat g_1 & \hat g_1 & \hat g_1 
		& b^1_{-1,-1,-1} & 8 \\ \hline
W_2^4 & g_1 & g_1 & g_2 & \hat g_1 & \hat g_1 & \hat g_2
		& b^1_{1,-1,-1}  & 8 \\ \hline
\end{array}$$

\medskip
\noindent We compute now the liftings.

\begin{algo} $V=V_2^4$.
We take $a_1$ and $a_2$ as in \ref{ta:lev}. Hence we have
\begin{align*}
\Dt(a_1)&=g_1\otimes a_1+a_1\otimes 1,
	\quad g_1a_1=-a_1g_1,\quad g_2a_1=a_1g_2; \\
\Dt(a_2)&=g_2\otimes a_2+a_2\otimes 1,
	\quad g_1a_2=-a_2g_1,\quad g_2a_2=-a_2g_2.
\end{align*}
Since $a_1^2\in\prim A$, we must have $a_1^2=0$ by \ref{rm:ebu}.
For the same reason, $a_2^2=0$. Take now
$$a_3=\Ad_{a_1}(a_2)=a_1a_2+g_1a_2\Ss(a_1)=a_1a_2+a_2a_1,$$
then
\begin{align*}
\Dt(a_3)&=g_1g_2\otimes a_1a_2+g_1a_2\otimes a_1+a_1g_2\otimes a_2
	+a_1a_2\otimes 1 \\
&=g_1g_2\otimes a_1a_2-a_2g_1\otimes a_1
	+a_1g_2\otimes a_2+a_1a_2\otimes 1.
\end{align*}
It is straightforward to see that $a_3^2\in\prim A$, and thus $a_3^2=0$.
Furthermore, it is clear from the definition of $a_3$ that
$$a_1a_3=a_3a_1,\ a_2a_3=a_3a_2,$$
and hence the unique lifting is the bosonization $\toba {V_2^2}\#\k\Gm$.
\end{algo}

\begin{algo} $V=W_2^1$.
As in \ref{ta:lev}, we can take $a_1$, $a_2$ and $a_3$ such that
$$\Dt(a_i)=g_1\otimes a_i+a_i\otimes 1,\quad
g_1a_i=-a_ig_1,\quad g_2a_i=a_ig_2.$$
We have $a_i^2\in\prim A$
and $a_ia_j+a_ja_i\in\prim A$, whence
$a_i^2=a_ia_j+a_ja_i=0$ $\forall i,j$ by \ref{rm:ebu}.
Hence, there is only one lifting:
the bosonization $\toba{W_2^1}\#\k\Gm$.
\end{algo}

\begin{algo} $V=W_2^2$.
This case is the same as the previous one, with the only difference
that $g_2a_3=-a_3g_2$, while $g_2a_i=a_ig_2$ for $i=1,2$. The same
arguments say that there is no other lifting than the bosonization
$\toba{W_2^2}\#\k\Gm$.
\end{algo}

\begin{algo} $V=W_2^3$.
We have now $g_1a_i=-a_ig_1$ and $g_2a_i=a_ig_2$. The difference with
the case $V=W_2^1$ is now that $\Dt(a_3)=g_1g_2\otimes a_3+a_3\otimes 1$.
As before, $a_i^2\in\prim A$ and hence $a_i^2=0$ by \ref{rm:ebu}.
Furthermore, $a_1a_2+a_2a_1\in\prim A$, from where $a_1a_2=-a_2a_1$.
However $a_ia_3+a_3a_i\in\primg A{g_2}$ for $i=1,2$ and hence, there
exist $\lmb_1,\ \lmb_2\in\k$ such that $a_ia_3+a_3a_i=\lmb_i(g_2-1)$.
Now, let $W'=(a_1,a_2)$ and let $T\in\GL_2(\k)$, $T=(T_{ij})$.
If we replace the basis $\{a_1,a_2\}$ by $\{Ta_1,Ta_2\}$ 
(we denote also by $T$ the automorphism of $W'$ with matrix $T$
in the basis $\{a_1,a_2\}$), then the constants
${\lmb_1\choose \lmb_2}$ are replaced by $T{\lmb_1\choose \lmb_2}$.
Hence, there are $2$ non isomorphic liftings: they can be
characterized by the rank of the morphism $W'\to\k(g_2-1)$,
$a\mapsto aa_3+a_3a$, which can be $0$ or $1$. Thus we have
\begin{enumerate}
\item The bosonization $\toba{W_2^3}\#\k\Gm$.
\item The algebra with the above relations plus
	$a_1a_3+a_3a_1=g_2-1$, 
	$a_2a_3+a_3a_2=0$.
\end{enumerate}
\end{algo}

\begin{algo} $V=W_2^4$.
We have here $\Dt(a_i)=g_1\otimes a_i+a_i\otimes 1$ for $i=1,2$, and
$\Dt(a_3)=g_2\otimes a_3+a_3\otimes 1$, plus $g_1a_i=-a_ig_1$,
$g_2a_i=a_ig_2$ for $i=1,2$, $g_1a_3=a_3g_1$, $g_2a_3=-a_3g_2$.
As before, $a_i^2=0$ for $i=1,2,3$, and also $a_1a_2+a_2a_1=0$.
Furthermore, $a_ia_3-a_3a_i\in\primg A{g_1g_2}$, from where
$a_ia_3-a_3a_1=\lmb_i(g_1g_2-1)$ for $i=1,2$.
These last conditions may add some unwanted relation which would
decrease the dimension of $A$. To see which $\lmb_i$'s give a
$32$-dimensional algebra, we use the diamond lemma:
\begin{align*}
g_1a_3a_1 &= g_1a_1a_3+\lmb_1g_1(g_1g_2-1) = -a_1a_3g_1+\lmb_1g_2-\lmb_1g_1 \\
g_1a_3a_1 &= -a_3a_1g_1 = -a_1a_3g_1-\lmb_1(g_1g_2-1)g_1 =
	-a_1a_3g_1-\lmb_1g_2+\lmb_1g_1,
\end{align*}
from where $\lmb_1=0$. An analogous work with the monomial
$g_1a_3a_2$ tells that $\lmb_2=0$, and there is only one lifting:
the bosonization $\toba{W_2^4}\#\k\Gm$.
\end{algo}

\subsection{$\Gm=C_4$} \label{sz4}
By similar considerations as for $\Gm=C_2\times C_2$, the rank must be $2$
or $3$. Here $\Aut(\Gm)\simeq C_2$, generated by $f(g)=g^3$. The list of
\ydm s giving Nichols algebras of dimension a power of $2$ is
$$\begin{array}{|c|l|l|c|c|} \hline
\multicolumn{5}{|c|}{\mbox{TABLE 5:
	modules $/\Aut(\Gm)$, $\Gm=C_4$}} \\ \hline
M(h,\chi) & h & \chi & (M(h,\chi))^f & \dim\toba{M(h,\chi)} \\ \hline
Y_3^1 & g & \hat g              & Y_3^8 & 4 \\ \hline
Y_3^2 & g & \hat g^2            & Y_3^7 & 2 \\ \hline
Y_3^3 & g & \hat g^3            & Y_3^6 & 4 \\ \hline
Y_3^4 & g^2 & \hat g            & Y_3^5 & 2 \\ \hline
Y_3^5 & g^2 & \hat g^3          & Y_3^4 & 2 \\ \hline
Y_3^6 & g^3 & \hat g            & Y_3^3 & 4 \\ \hline
Y_3^7 & g^3 & \hat g^2          & Y_3^2 & 2 \\ \hline
Y_3^8 & g^3 & \hat g^3          & Y_3^1 & 4 \\ \hline
\end{array}$$

Hence, for rank $2$ we have the following non isomorphic
possibilities (we exclude the cases $Y_3^i\oplus Y_3^j$ with
$i,j\in\{1,3,6,8\}$ by virtue of \ref{ex:qls}).

$$\begin{array}{|c|l|l|l|l|r|r|r|r|l|c|} \hline
\multicolumn{11}{|c|}{\mbox{TABLE 6:
	Rank $2$, $\Gm=C_4$}} \\ \hline
V     & h_1 & h_2 & \chi_1 & \chi_2 & b_{11} & b_{12} &
	b_{21} & b_{22} & b_{ij} & \dim\toba V \\ \hline
V_3^1 & g & g & \hat g & \hat g^2
		&  \rmu & -1 &  \rmu & -1 & [b^4_-]^{\tau} & 64 \\ \hline
V_3^2 & g & g & \hat g^3 & \hat g^2
		& -\rmu & -1 & -\rmu & -1 & [b^4_+]^{\tau} & 64 \\ \hline
V_3^3 & g & g & \hat g^2 & \hat g^2
		& -1 & -1 & -1 & -1 & b^1_{-1}& 4 \\ \hline
V_3^4 & g & g^3 & \hat g^2 & \hat g^3
		& -1 & -\rmu & -1 &  \rmu & [b^5_-]^{\rho\tau}
		& \ge 24 \\ \hline
V_3^5 & g & g^3 & \hat g^2 & \hat g
		& -1 &  \rmu & -1 & -\rmu & [b^5_+]^{\rho\tau}
		& \ge 24 \\ \hline
V_3^6 & g & g^3 & \hat g^2 & \hat g^2
		& -1 & -1 & -1 & -1 & b^1_{-1} & 4 \\ \hline
V_3^7 & g & g^2 & \hat g & \hat g
		&  \rmu &  \rmu & -1 & -1 & b^4_+ & 64 \\ \hline
V_3^8 & g & g^2 & \hat g & \hat g^3
		&  \rmu & -\rmu & -1 & -1 & b^5_+ & \ge 24 \\ \hline
V_3^9 & g & g^2 & \hat g^3 & \hat g
		& -\rmu &  \rmu & -1 & -1 & b^5_- & \ge 24 \\ \hline
V_3^{10}& g & g^2 & \hat g^3 & \hat g^3
		& -\rmu & -\rmu & -1 & -1 & b^4_- & 64 \\ \hline
V_3^{11}& g & g^2 & \hat g^2 & \hat g
		& -1 &  \rmu &  1 & -1 & [b^3_{++}]^{\rho} & 16 \\ \hline
V_3^{12}& g & g^2 & \hat g^2 & \hat g^3
		& -1 & -\rmu &  1 & -1 & [b^3_{-+}]^{\rho} & 16 \\ \hline
V_3^{13}& g^2 & g^2 & \hat g & \hat g
		& -1 & -1 & -1 & -1 & b^1_{-1} & 4 \\ \hline
V_3^{14}& g^2 & g^2 & \hat g & \hat g^3
		& -1 & -1 & -1 & -1 & b^1_{-1} & 4 \\ \hline
\end{array}$$

\medskip
Thus, we have no rank $2$ modules over $C_4$ giving $8$-dimensional
Nichols algebras.
For rank $3$ the possible subcases are $V_3^3$, $V_3^6$, $V_3^{13}$ and
$V_3^{14}$. We have hence the following cases:

$$\begin{array}{|c|l|l|l|l|l|l|c|c|} \hline
\multicolumn{9}{|c|}{\mbox{TABLE 7:
	Rank $3$, $\Gm=C_4$}} \\ \hline
V     & h_1 & h_2 & h_3 & \chi_1 & \chi_2 & \chi_3 
	& b_{ij} & \dim\toba V \\ \hline
W_3^1 & g & g & g & \hat g^2 & \hat g^2 & \hat g^2 & b^1_{-1,-1,-1} & 8 \\ \hline
W_3^2 & g & g & g^3 & \hat g^2 & \hat g^2 & \hat g^2 & b^1_{-1,-1,-1} & 8 \\ \hline
W_3^3 & g^2 & g^2 & g^2 & \hat g & \hat g & \hat g & b^1_{-1,-1,-1} & 8 \\ \hline
W_3^4 & g^2 & g^2 & g^2 & \hat g & \hat g & \hat g^3 & b^1_{-1,-1,-1} & 8 \\ \hline
\end{array}$$

\medskip
We compute now the liftings:

\begin{algo} $V=W_3^1$.
We take $a_i,\ i=1,2,3$ as in \ref{ta:lev}, i.e.
$\Dt(a_i)=g\otimes a_i+a_i\otimes 1$, $ga_i=-a_ig$.
Hence, by \ref{rm:ebu}, $a_i^2\in\primg A{g^2}$, and
$a_ia_j+a_ja_i\in\primg A{g^2}$, from where there exist
$\lmb_i, 1\le i\le 6$ such that
\begin{align*}
& a_1^2=\lmb_1(g^2-1),\  a_2^2=\lmb_2(g^2-1),\  a_3^2=\lmb_3(g^2-1), \\
& a_1a_2+a_2a_1=\lmb_4(g^2-1), \\
& a_1a_3+a_3a_1=\lmb_5(g^2-1),\\
& a_2a_3+a_3a_2=\lmb_6(g^2-1).
\end{align*}
Thus, we have a quadratic form $f:W=(a_1,a_2,a_3)\to\k(g^2-1)$,
$f(a)=a^2$. Furthermore, the basis $\{a_1,a_2,a_3\}$ of $W$ can
be replaced by any other basis of $W$, as a result of what we must
consider the equivalence class of the quadratic form $f$ in order
to have non isomorphic liftings. Since $\car\k\neq 2$ and $\k$ is
algebraically closed, $f$ is characterized by its rank,
and hence there are $4$ liftings (up to isomorphism):
\begin{enumerate}
\item $a_ia_j+a_ja_i=0$ for $i\neq j$, $a_i^2=0$ for $i=1,2,3$. This is
	the bosonization $\toba{W_3^1}\#\k\Gm$.
\item $a_ia_j+a_ja_i=0$ for $i\neq j$, $a_i^2=0$ for $i=1,2$, $a_3^2=(g^2-1)$.
\item $a_ia_j+a_ja_i=0$ for $i\neq j$, $a_1^2=0$, $a_i^2=(g^2-1)$ for $i=2,3$.
\item $a_ia_j+a_ja_i=0$ for $i\neq j$, $a_i^2=(g^2-1)$ for $i=1,2,3$.
\end{enumerate}
\end{algo}

\begin{algo} $V=W_3^2$.
Similarly, we have $a_1,a_2\in\primg Ag$, $a_3\in\primg A{g^3}$ and such
that $ga_i=-a_ig$ for $i=1,2,3$. Now, $a_ia_3+a_3a_i\in\prim A$ ($i=1,2$),
whence $a_ia_3+a_3a_i=0$ for $i=1,2$. However, there exist
$\lmb_i,\ 1\le i\le 4$ such that
$$a_i^2=\lmb_i(g^2-1)\ (i=1,2,3),\qquad a_1a_2+a_2a_1=\lmb_4(g^2-1).$$
Hence we have two quadratic forms:
\begin{align*}
f_1&:W_1=(a_1,a_2)\to\k(g^2-1),\ f_1(a)=a^2, \\
f_2&:W_2=(a_3)\to\k(g^2-1),\ f_2(a)=a^2.
\end{align*}
As before, we can replace $\{a_1,a_2\}$ by another basis of $W_1$,
and the same for $\{a_3\}$. As with $W_3^1$, $f_1$ and $f_2$ are
characterized by their ranks, whence there are $3$ possibilities
for $f_1$ and $2$ for $f_2$, and finally there are $6$ liftings:
\begin{enumerate}
\item $a_ia_j+a_ja_i=0$ for $i\neq j$, $a_1^2=0$, $a_2^2=0$, $a_3^2=0$.
	This is the bosonization $\toba{W_3^2}\#\k\Gm$.
\item $a_ia_j+a_ja_i=0$ for $i\neq j$, $a_1^2=0$, $a_2^2=(g^2-1)$, $a_3^2=0$.
\item $a_ia_j+a_ja_i=0$ for $i\neq j$, $a_1^2=(g^2-1)$, $a_2^2=(g^2-1)$,
	$a_3^2=0$.
\item $a_ia_j+a_ja_i=0$ for $i\neq j$, $a_1^2=0$, $a_2^2=0$, $a_3^2=(g^2-1)$.
\item $a_ia_j+a_ja_i=0$ for $i\neq j$, $a_1^2=0$, $a_2^2=(g^2-1)$,
	$a_3^2=(g^2-1)$.
\item $a_ia_j+a_ja_i=0$ for $i\neq j$, $a_1^2=(g^2-1)$, $a_2^2=(g^2-1)$,
	$a_3^2=(g^2-1)$.
\end{enumerate}
\end{algo}

\begin{algo} $V=W_3^3$.
The situation here is easier: we have $a_i\in\primg A{g^2}$ and
$ga_i=\rmu a_ig$ for $i=1,2,3$. Since $a_i^2\in\prim A$ and
$a_ia_j+a_ja_i\in\prim A$ for all $i,j$, by \ref{rm:ebu} we have
$a_i^2=a_ia_j+a_ja_i=0$, and there is only one
lifting: the bosonization $\toba{W_3^3}\#\k\Gm$.
\end{algo}

\begin{algo} $V=W_3^4$.
Similar to the previous case, the only difference being that
$ga_i=\rmu a_ig$ for $i=1,2$ and $ga_3=-\rmu a_3g$. There is only one lifting:
the bosonization $\toba{W_3^3}\#\k\Gm$.
\end{algo}

\bigskip

\section{$\Gm$ of order $8$} \label{ss:or8}
By \ref{lm:nht} we have $\dim V\le 2$. If $\dim V=1$, we have
$V=M(h,\chi)$, $h\in Z(\Gm)$ and $\chi$ a character such that
$\chi(h)$ has order $4$. If $\Gm$ is non abelian, we have
$Z(\Gm)\subseteq [\Gm;\Gm]$, whence any irreducible representation
of degree $1$ vanishes on the center, and $\chi(h)=1$. This
tells that there are no $32$-dimensional pointed Hopf algebras with
non abelian coradical of dimension $8$.
\begin{rem}\label{rm:cdna}
In turn this explains why there are no $16$-dimensional pointed Hopf
algebras with non abelian coradical (see \cite{cdr}).
\end{rem}

If $\dim V=2$. Then, by
\cite[Prop. 3.1.11]{ag}, $V$ comes from the abelian case, and hence
it has a matrix $(b_{ij})$. By \ref{ex:qls}, $V$ must be a QLS over
some subgroup $\Gm'\subseteq\Gm$, with $b_{11}=b_{22}=-1$.
If $\Gm$ is non abelian, there are in principle three possibilities
for $V$ to be $2$-dimensional:
\begin{algo}\label{al:casos}
Decomposition of $V$, $\Gm$ non abelian.
\begin{enumerate}
\item $V=M(h_1,\chi_1)\oplus M(h_2,\chi_2)$, where the $h_i$'s are central
	in $\Gm$ and the $\chi_i$'s are characters, \label{ca:d21}
\item $V=M(h,\chi)$, where $[\Gm:\Gm_h]=2$ and $\chi$ is a character of
	$\Gm_h$, \label{ca:d22}
\item $V=M(h,\rho)$, where $h$ is central in $\Gm$ and $\dim\rho=2$.
	\label{ca:d23}
\end{enumerate}
\end{algo}
Case \ref{ca:d21} does not arise for the same reasons as above: any
irreducible representation of degree $1$ of a non abelian group of
order $8$ vanishes on the center, and we would have $b_{11}=b_{22}=1$,
contradicting \ref{rm:qn1}. We shall proceed as follows:
\begin{tactic}\label{ta:nal1}
In case \ref{ca:d22}, let $x$ be a generator of the space affording $\chi$.
Note that $\chi(h)=b_{11}=-1$. Let $t\not\in\Gm_h$.
The conjugacy class of $h$ is
$\{h,\alai ht=tht^{-1}\}$, and since
$[\Gm:\Gm_h]=2$ it is easy to prove that $\alai ht\in\Gm_h$.
Notice that $t^{-1}\not\in\Gm_h$, whence $t^{-1}ht=tht^{-1}=\alai ht$.
Let $y=t\ganda x$ and let $\zeta=\chi(\alai ht)$.
Then $\dt(y)=\alai ht\otimes y$, and $c$ can be computed as follows:
\begin{align*}
c(x\otimes x)&=(h\ganda x)\otimes x=-x\otimes x, \\
c(x\otimes y)&=(h\ganda y)\otimes x=(ht\ganda x)\otimes x
	=(tt^{-1}ht\ganda x)\otimes x \\
	&\qquad\qquad =(\zeta t\ganda x)\otimes x
	=\zeta y\otimes x, \\
c(y\otimes x)&=(\alai ht\ganda x)\otimes y=\zeta x\otimes y, \\
c(y\otimes y)&=(\alai htt\ganda x)\otimes y=(th\ganda x)\otimes y=
	-y\otimes y.
\end{align*}
Therefore, the matrix of $c$ is $\mdpd {-1}{\zeta}{\zeta}{-1}$, whence
$1=b_{12}b_{21}=\zeta^2$, i.e. $\zeta=\pm 1$.
To compute the liftings, let $a,b$ be liftings of $x,y$ as in
\ref{ta:lev}. Then $hah^{-1}=-a$ and $b=tat^{-1}$. The action
of $\Gm$ on $A_1$ may be computed from the action of $\Gm$
on $V$, as explained in \ref{ta:lev}.
\end{tactic}

\begin{tactic} \label{ta:nal2}
In case \ref{ca:d23}, let $x$ be a generator of the space affording
$\rho$, and let $t\in\Gm$ such that $y=\rho(t)(x)$ is linearly
independent with $x$.
Since $\rho$ is irreducible and $h$ is central, $\rho(h)=\zeta\id$
for some $\zeta\in\unid{\k}$. Then the matrix of $c$ is
$\mdpd {\zeta}{\zeta}{\zeta}{\zeta}$,
from where $\zeta=b_{11}=-1$. Let $a,b$ be liftings of $x,y$
as in \ref{ta:lev}. Then $b=tat^{-1}$ and $hbh^{-1}=-b$.
\end{tactic}

\begin{rem}
These cases do not arise when classifying pointed Hopf algebras of
dimension $p^5$, $p$ an odd prime number (see \cite{gp5}), and there
are no pointed Hopf algebras of dimension $p^5$ with non commutative
coradical of dimension $p^3$.
\end{rem}

We consider now the case $\Gm$ abelian.
\begin{algo}\label{al:uodl}
Use of diamond lemma

If $\Gm$ is abelian and $\dim V=1$, then $V=M(h,\chi)$
and, as in \ref{ta:lev}, we have $a=a_1\in\primg Ah$.
Furthermore, for $\toba V$ to be $4$-dimensional,
we must have $\chi(h)$ of order $4$.
By \ref{rm:ebu}, $a^4\in\primg A{h^4}=\k(h^4-1)$ ($\chi(h)$ having
order $4$, this implies that $h^4$ does not coincide neither
with $h$, $h^2$ nor $h^3$).
We have hence $A$ generated by $\Gm$ and $a$, with relations
$$tat^{-1}=\chi(t)a\ \forall t\in\Gm,\qquad
	a^4=\lmb(h^4-1).$$
If $h^4\neq 1$, taking a suitable scalar multiple of $a$ we
may suppose that $\lmb\in\ceun$. We use the diamond lemma
to prove that $A$ is $32$-dimensional, but this adds some
condition on $\lmb$. Let $t\in\Gm$. We have
\begin{align*}
ta^4&=\chi^4(t)a^4t=\chi^4(t)\lmb(h^4-1)t, \\
ta^4&=t\lmb(h^4-1)=\lmb(h^4-1)t,
\end{align*}
whence
\begin{equation}\label{eq:cd1}
\lmb(\chi^4(t)-1)(h^4-1)=0\ \forall t\in\Gm.
\end{equation}
If $h^4=1$ we can take $\lmb=0$, and hence condition \eqref{eq:cd1}
can be written in the more compact form
\begin{equation}\label{eq:cd2}
\lmb(\chi^4-1)=0.
\end{equation}

If $\Gm$ is abelian and $\dim V=2$, then
$V=M(h_1,\chi_1)\oplus M(h_2,\chi_2)$. We have
as in \ref{ta:lev} $a_i\in\primg A{h_i}$, $ta_it^{-1}=\chi_i(t)a_i$.
As in \ref{rm:ebu}, $a_i^2\in\primg A{h_i^2}$ and
$\Ad_{a_1}(a_2)=a_1a_2-\chi_2(h_1)a_2a_1\in\primg A{h_1h_2}$.
Since $\toba V$ is $4$-dimensional, by \ref{ex:qls} it is a QP,
and we have $\chi_2(h_1)\chi_1(h_2)=1$, $\chi_i(h_i)=-1$. Now,
$h_1h_2\neq h_1$ and $h_1h_2\neq h_2$, for if not $h_2=1$ or
$h_1=1$, contradicting $\chi_i(h_i)=-1$. Moreover, if $h_1^2=h_2$,
we would have
$$-1=\chi_2(h_2)=\chi_2(h_1^2)=(\chi_2(h_1))^2=(\chi_1(h_2))^{-2}
	=(\chi_1(h_1))^{-4}=(-1)^{-4}=1,$$
which is absurd. Thus $h_1^2\neq h_2$, and analogously $h_2^2\neq h_1$.

\noindent Hence, $\primg A{h_i^2}=\k(h_i^2-1)$ and
$\primg A{h_1h_2}=\k(h_1h_2-1)$, and we have
\begin{equation}\label{eq:edl}
a_i^2=\lmb_i(h_i^2-1)\ (i=1,2),\qquad
a_1a_2-\chi_2(h_1)a_2a_1=\lmb_3(h_1h_2-1).
\end{equation}
Similarly to \eqref{eq:cd2}, we get
\begin{equation}\label{eq:cd3}
\lmb_1(\chi_1^2-1)=0,\quad \lmb_2(\chi_2^2-1)=0.
\end{equation}
Doing the same with the monomials $ta_2a_1$, $t\in\Gm$, we
arrive to the condition
\begin{equation}\label{eq:cd4}
\lmb_3(\chi_1\chi_2-1)=0.
\end{equation}
\end{algo}

\medskip
We consider now group by group.
\subsection{$\Gm=C_2\times C_2\times C_2$} \label{sz2z2z2}
There are no $1$-dimensional modules $M(h,\chi)$ in $\ydskg$
s.t. $\chi(h)$ has order $4$, whence the rank must be $2$.
Now, the \ydm s giving $2$-dimensional Nichols algebras are those
$M(h,\chi)$ such that $\chi(h)=-1$, but
$\Aut(\Gm)\simeq \GL_3(\ZZ/2)$, and this group acts transitively
on these modules. Hence we can suppose that $M_1=M(g_1,\hat g_1)$
(we are using the notation at the end of \ref{ss:bdac}).
The subgroup of $\Aut(\Gm)$ which fixes $M_1$ is isomorphic to
$\SS_3$, and is generated by
$$f_1=(g_1\mapsto g_1,g_2\mapsto g_2g_3,g_3\mapsto g_3),\quad
f_2=(g_1\mapsto g_1,g_2\mapsto g_2,g_3\mapsto g_2g_3).$$
Thus we consider the action of this subgroup on the modules and
choose one element in each orbit:
\begin{align*}
& Y_4^1=M(g_1,\hat g_1),\ Y_4^2=M(g_1,\hat g_1\hat g_2),
\ Y_4^3=M(g_1g_2,\hat g_1),\ Y_4^4=M(g_1g_2,\hat g_1\hat g_3), \\
& Y_4^5=M(g_1g_2,\hat g_2),\ Y_4^6=M(g_2,\hat g_2),
\ Y_4^7=M(g_2,\hat g_1\hat g_2).
\end{align*}
In order to give non isomorphic modules after taking the sums
$Y_4^1\oplus Y_4^i$, we must take into account the possibility
of switching the order, as in lemma \ref{lm:iso}. The modules
$Y_4^1\oplus Y_4^5$ and $Y_4^1\oplus Y_4^7$ give isomorphic
algebras when its Nichols algebras are bosonized by means of
the automorphism
$$(g_1\mapsto g_2,\ g_2\mapsto g_1g_2,\ g_3\mapsto g_3),$$
and we consider only $Y_4^1\oplus Y_4^5$. Hence, the list is
$$\begin{array}{|l|l|l|l|l|r|r|r|r|l|c|} \hline
\multicolumn{11}{|c|}{\mbox{TABLE 8:
	Rank $2$, $\Gm=C_2\times C_2\times C_2$}} \\ \hline
V     & h_1 & h_2 & \chi_1 & \chi_2 & b_{11} & b_{12} & b_{21} 
	& b_{22} & b_{ij} & \dim\toba V \\ \hline
V_4^1 & g_1 & g_1 & \hat g_1 & \hat g_1  & -1&-1&-1&-1          
		& b^1_{-1} & 4 \\ \hline
V_4^2 & g_1 & g_1 & \hat g_1 & \hat g_1\hat g_2 & -1&-1&-1&-1   
		& b^1_{-1} & 4 \\ \hline
V_4^3 & g_1 & g_1g_2 & \hat g_1 & \hat g_1 & -1&-1&-1&-1        
		& b^1_{-1} & 4 \\ \hline
V_4^4 & g_1 & g_1g_2 & \hat g_1 & \hat g_1\hat g_3 & -1&-1&-1&-1
		& b^1_{-1} & 4 \\ \hline
V_4^5 & g_1 & g_1g_2 & \hat g_1 & \hat g_2 & -1&1 &-1&-1        
		& [b^2]^{\rho} & 8 \\ \hline
V_4^6 & g_1 & g_2 & \hat g_1 & \hat g_2 & -1&1 &1 &-1           
		& b^1_{1} & 4 \\ \hline
\end{array}$$

We compute the liftings:

\begin{algo} $V=V_4^1$.
We have, as in \ref{ta:lev}, $a_1,a_2\in\primg A{g_1}$, and
$g_1a_i=-a_ig_1$, $g_2a_i=a_ig_2$, $g_3a_i=a_ig_3$ for $i=1,2$.
By \ref{rm:ebu}, we have $a_i^2\in\prim A$ and also
$a_1a_2+a_2a_1\in\prim A$,
whence $a_i^2=a_1a_2+a_2a_1=0$, and there is only one lifting: the
bosonization $\toba{V_4^1}\#\k\Gm$.
\end{algo}

\begin{algo} $V=V_4^2$.
It is similar to the previous one, the only difference being that
$g_2a_2=-a_2g_2$. The unique lifting is
the bosonization $\toba{V_4^2}\#\k\Gm$
\end{algo}

\begin{algo} $V=V_4^3$.
Here $a_1\in\primg A{g_1}$, $a_2\in\primg A{g_1g_2}$
and $\Gm$ acts by 
$g_1a_i=-a_ig_1$, $g_2a_i=a_ig_2$, $g_3a_i=a_ig_3$. As before,
$a_i^2\in\prim A$, whence $a_i^2=0$ for $i=1,2$. However,
$a_1a_2+a_2a_1\in\primg A{g_2}$, from where as in \eqref{eq:edl}
we have
$a_1a_2+a_2a_1=\lmb_3(g_2-1)$.
Taking a suitable scalar multiple of $a_1$, we may suppose that
$\lmb_3=0$ or $\lmb_3=1$. There are hence two liftings:
\begin{enumerate}
\item $a_i^2=0$, $a_1a_2+a_2a_1=0$. This is the bosonization
	$\toba{V_4^3}\#\k\Gm$.
\item $a_i^2=0$, $a_1a_2+a_2a_1=(g_2-1)$.
\end{enumerate}
\end{algo}

\begin{algo} $V=V_4^4$.
It is similar to $V_4^3$. We have $a_1\in\primg A{g_1}$,
$a_2\in\primg A{g_1g_2}$ and $\Gm$ acts now by
$g_1a_i=-a_ig_1$, $g_2a_i=a_ig_2$, $g_3a_1=a_1g_3$, $g_3a_2=-a_2g_3$.
As before, $a_i^2=0$, $a_1a_2+a_2a_1=\lmb_3(g_2-1)$.
Now, \eqref{eq:cd4} tells that $\lmb_3=0$. There is only one
lifting: the bosonization $\toba{V_4^4}\#\k\Gm$.
\end{algo}

\begin{algo} $V=V_4^6$.
Now $a_1\in\primg A{g_1}$, $a_2\in\primg A{g_2}$ and $\Gm$ acts by
$$g_1a_1=-a_1g_1,\quad g_1a_2=a_2g_1,\quad g_2a_1=a_1g_2,\quad
g_2a_2=-a_2g_2,\quad g_3a_i=a_ig_3.$$
As before, $a_i^2\in\prim A$, whence $a_i^2=0$.
Furthermore, $a_1a_2-a_2a_1\in\primg A{g_1g_2}$, from where
$a_1a_2-a_2a_1=\lmb_3(g_1g_2-1)$.
Now, \eqref{eq:cd4} tells that $\lmb_3=0$. There is only one
lifting: the bosonization $\toba{V_4^6}\#\k\Gm$.
\end{algo}

\subsection{$\Gm=C_2\times C_4$} \label{sz2z4}
We have, using the notation at the end of \ref{ss:bdac},
$\Aut(\Gm)\simeq\DD_4$ with generators
$$f_1=(g_1\mapsto g_1g_2^2,\ g_2\mapsto g_1g_2),\quad
	f_2=(g_1\mapsto g_1g_2^2,\ g_2\mapsto g_2),$$
with relations $f_1^4=f_2^2=1$, $f_1f_2=f_2f_1^3$. Hence there are $3$
orbits in $\Gm$ under the action of $\Aut(\Gm)$:
$\{g_2,g_1g_2,g_2^3,g_1g_2^3\}$, $\{g_1,g_1g_2^2\}$ and $\{g_2^2\}$.
There are $7$ orbits of irreducible \ydm s giving Nichols algebras
of dimension a power of $2$. The following is a list of one element
per orbit and the subgroup of $\Aut(\Gm)$ fixing the chosen element.
$$\begin{array}{|c|l|l|c|c|} \hline
\multicolumn{5}{|c|}{\mbox{TABLE 9:
	modules$/\Aut(\Gm)$, $\Gm=C_2\times C_4$}} \\ \hline
M(h,\chi) & h & \chi & \dim\toba M & (\Aut(\Gm))_{M} \\ \hline
Y_5^1 & g_2 & \hat g_2                  & 4 & 1 \\ \hline
Y_5^2 & g_2 & \hat g_2^2                & 2 & (f_2) \\ \hline
Y_5^3 & g_2 & \hat g_2^3                & 4 & 1 \\ \hline
Y_5^4 & g_2 & \hat g_1\hat g_2^2        & 2 & (f_2) \\ \hline
Y_5^5 & g_1 & \hat g_1\hat g_2^3        & 2 & (f_2f_1) \\ \hline
Y_5^6 & g_1 & \hat g_1                  & 2 & (f_1^2) \\ \hline
Y_5^7 & g_2^2 & \hat g_2                & 2 & (f_1f_2) \\ \hline
\end{array}$$
Hence there are two rank $1$ Nichols algebras of dimension $4$ over
$C_2\times C_4$. We go now to rank $2$, i.e. $V=M_1\oplus M_2$.
Since $\dim\toba{M_i}=2$, we can suppose that $M_1$ is one of
$Y_5^2$, $Y_5^4$, $Y_5^5$, $Y_5^6$ or $Y_5^7$.
We give the list, which is the hardest to compute (we take
into account \ref{lm:iso}).
$$\begin{array}{|c|l|l|l|l|r|r|r|r|l|c|} \hline
\multicolumn{11}{|c|}{\mbox{TABLE 10:
	Rank $2$, $\Gm=C_2\times C_4$}} \\ \hline
V     & h_1 & h_2 & \chi_1 & \chi_2 
	& b_{11} & b_{12} & b_{21} & b_{22} & b_{ij} & \dim\toba V \\ \hline
V_5^1 & g_2 & g_2 & \hat g_2^2 & \hat g_2^2
	& -1 & -1 & -1 & -1 & b^1_{-1} & 4 \\ \hline
V_5^2 & g_2 & g_2 & \hat g_2^2 & \hat g_1\hat g_2^2
	& -1 & -1 & -1 & -1 & b^1_{-1} & 4 \\ \hline
V_5^3 & g_2 & g_2 & \hat g_1\hat g_2^2 & \hat g_1\hat g_2^2
	& -1 & -1 & -1 & -1 & b^1_{-1} & 4 \\ \hline
V_5^4 & g_2 & g_1g_2 & \hat g_2^2 & \hat g_2^2
	& -1 & -1 & -1 & -1 & b^1_{-1} & 4 \\ \hline
V_5^5 & g_2 & g_1g_2 & \hat g_2^2 & \hat g_1
	& -1 &  1 & -1 & -1 & [b^2]^{\rho} & 8 \\ \hline
V_5^6 & g_2 & g_1g_2 & \hat g_1\hat g_2^2 & \hat g_1
	& -1 &  1 &  1 & -1 & b^1_{1} & 4 \\ \hline
V_5^7 & g_2 & g_2^3 & \hat g_2^2 & \hat g_2^2
	& -1 & -1 & -1 & -1 & b^1_{-1} & 4 \\ \hline
V_5^8 & g_2 & g_2^3 & \hat g_1\hat g_2^2 & \hat g_2^2
	& -1 & -1 & -1 & -1 & b^1_{-1} & 4 \\ \hline
V_5^{9} & g_2 & g_2^3 & \hat g_1\hat g_2^2 & \hat g_1\hat g_2^2
	& -1 & -1 & -1 & -1 & b^1_{-1} & 4 \\ \hline
V_5^{10} & g_2 & g_1 & \hat g_2^2 & \hat g_1
	& -1 &  1 &  1 & -1 & b^1_{1} & 4 \\ \hline
V_5^{11} & g_2 & g_1 & \hat g_2^2 & \hat g_1\hat g_2^2
	& -1 & -1 &  1 & -1 & b^2 & 8 \\ \hline
V_5^{12} & g_2 & g_1 & \hat g_1\hat g_2^2 & \hat g_1
	& -1 &  1 & -1 & -1 & [b^2]^{\rho} & 8 \\ \hline
V_5^{13} & g_2 & g_1 & \hat g_1\hat g_2^2 & \hat g_1\hat g_2^2
	& -1 & -1 & -1 & -1 & b^1_{-1} & 4 \\ \hline
V_5^{14} & g_2 & g_1 & \hat g_2^2 & \hat g_1\hat g_2^3
	& -1 & -\rmu &  1 & -1 & [b^3_{-+}]^{\rho} & 16 \\ \hline
V_5^{15} & g_2 & g_1 & \hat g_2^2 & \hat g_1\hat g_2
	& -1 &  \rmu &  1 & -1 & [b^3_{++}]^{\rho} & 16 \\ \hline
V_5^{16} & g_2 & g_1 & \hat g_1\hat g_2^2 & \hat g_1\hat g_2^3
	& -1 & -\rmu & -1 & -1 & [b^3_{--}]^{\rho} & 16 \\ \hline
V_5^{17} & g_2 & g_1 & \hat g_1\hat g_2^2 & \hat g_1\hat g_2
	& -1 &  \rmu & -1 & -1 & [b^3_{+-}]^{\rho} & 16 \\ \hline
V_5^{18} & g_2 & g_2^2 & \hat g_2^2 & \hat g_1\hat g_2
	& -1 &  \rmu &  1 & -1 & [b^3_{++}]^{\rho} & 16 \\ \hline
V_5^{19} & g_2 & g_2^2 & \hat g_2^2 & \hat g_1\hat g_2^3
	& -1 & -\rmu &  1 & -1 & [b^3_{-+}]^{\rho} & 16 \\ \hline
V_5^{20} & g_2 & g_2^2 & \hat g_1\hat g_2^2 & \hat g_2
	& -1 &  \rmu &  1 & -1 & [b^3_{++}]^{\rho} & 16 \\ \hline
V_5^{21} & g_2 & g_2^2 & \hat g_1\hat g_2^2 & \hat g_2^3
	& -1 & -\rmu &  1 & -1 & [b^3_{-+}]^{\rho} & 16 \\ \hline
V_5^{22} & g_2^2 & g_2^2 & \hat g_2 & \hat g_2
	& -1 & -1 & -1 & -1 & b^1_{-1} & 4 \\ \hline
V_5^{23} & g_2^2 & g_2^2 & \hat g_2 & \hat g_1\hat g_2
	& -1 & -1 & -1 & -1 & b^1_{-1} & 4 \\ \hline
V_5^{24} & g_2^2 & g_2^2 & \hat g_2 & \hat g_2^3
	& -1 & -1 & -1 & -1 & b^1_{-1} & 4 \\ \hline
V_5^{25} & g_2^2 & g_1 & \hat g_2 & \hat g_1
	& -1 &  1 &  1 & -1 & b^1_{1} & 4 \\ \hline
V_5^{26} & g_2^2 & g_1 & \hat g_2 & \hat g_1\hat g_2^3
	& -1 & -1 &  1 & -1 & b^2 & 8 \\ \hline
V_5^{27} & g_2^2 & g_1 & \hat g_1\hat g_2^3 & \hat g_1\hat g_2^3
	& -1 & -1 & -1 & -1 & b^1_{-1} & 4 \\ \hline
V_5^{28} & g_2^2 & g_1 & \hat g_1\hat g_2^3 & \hat g_1\hat g_2
	& -1 & -1 & -1 & -1 & b^1_{-1} & 4 \\ \hline
V_5^{29} & g_2^2 & g_1 & \hat g_1\hat g_2 & \hat g_1
	& -1 &  1 & -1 & -1 & [b^2]^{\rho} & 8 \\ \hline
V_5^{30} & g_2^2 & g_1 & \hat g_2^3 & \hat g_1
	& -1 &  1 &  1 & -1 & b^1_1 & 4 \\ \hline
V_5^{31} & g_2^2 & g_1g_2^2 & \hat g_2 & \hat g_1
	& -1 &  1 & -1 & -1 & [b^2]^{\rho} & 8 \\ \hline
V_5^{32} & g_1 & g_1 & \hat g_1\hat g_2^3 & \hat g_1\hat g_2^3
	& -1 & -1 & -1 & -1 & b^1_{-1} & 4 \\ \hline
\multicolumn{11}{|c|}{\mbox{Continues in next page}} \\ \hline
\end{array}$$
$$\begin{array}{|c|l|l|l|l|r|r|r|r|l|c|} \hline
V       & h_1 & h_2 & \chi_1 & \chi_2 
	& b_{11} & b_{12} & b_{21} & b_{22} & b_{ij} & \dim\toba V \\ \hline
V_5^{33} & g_1 & g_1 & \hat g_1\hat g_2^3 & \hat g_1\hat g_2
	& -1 & -1 & -1 & -1 & b^1_{-1} & 4 \\ \hline
V_5^{34} & g_1 & g_1 & \hat g_1\hat g_2^3 & \hat g_1
	& -1 & -1 & -1 & -1 & b^1_{-1} & 4 \\ \hline
V_5^{35} & g_1 & g_1 & \hat g_1 & \hat g_1
	& -1 & -1 & -1 & -1 & b^1_{-1} & 4 \\ \hline
V_5^{36} & g_1 & g_1 & \hat g_1 & \hat g_1\hat g_2^2
	& -1 & -1 & -1 & -1 & b^1_{-1} & 4 \\ \hline
V_5^{37} & g_1 & g_1g_2^2 & \hat g_1\hat g_2^3 & \hat g_2
	& -1 &  1 &  1 & -1 & b^1_{1} & 4 \\ \hline
V_5^{38} & g_1 & g_1g_2^2 & \hat g_1\hat g_2^3 & \hat g_1
	& -1 & -1 &  1 & -1 & b^2 & 8 \\ \hline
V_5^{39} & g_1 & g_1g_2^2 & \hat g_1 & \hat g_1
	& -1 & -1 & -1 & -1 & b^1_{-1} & 4 \\ \hline
V_5^{40} & g_1 & g_1g_2^2 & \hat g_1 & \hat g_1\hat g_2^2
	& -1 & -1 & -1 & -1 & b^1_{-1} & 4 \\ \hline
\end{array}$$

\bigskip
We compute the liftings. For rank $1$ we have
\begin{algo} $V=Y_5^1$.
As in \ref{ta:lev}, we have $a\in\primg A{g_2}$, $g_1a=ag_1$,
$g_2a=\rmu ag_2$. By \ref{rm:ebu}, we have $a^4=0$ and there is
only one lifting: the bosonization $\toba{Y_5^1}\#\k\Gm$.
\end{algo}

\begin{algo} $V=Y_5^3$.
Analogous to the previous one, the only difference being
that $g_2a=-\rmu ag_2$. There is only one lifting: the
bosonization $\toba{Y_5^3}\#\k\Gm$.
\end{algo}

\medskip
We go now to rank $2$. As explained in \ref{al:uodl},
we have $a_1,a_2$ with 
$a_i\in\primg A{h_i}$, $ta_it^{-1}=\chi_i(t)a_i$,
$$a_1^2=\lmb_1(h_1^2-1),\quad a_2^2=\lmb_2(h_2^2-1),\quad
a_1a_2-\chi_2(h_1)a_2a_1=\lmb_3(h_1h_2-1).$$
Before computing the liftings, we notice some facts.
\begin{rem}\label{rm:itr}
In certain cases, some element in $\Aut(\Gm)$ interchanges $M_1$ and
$M_2$. In this situation the lifting corresponding to $(\lmb_1,\lmb_2,\lmb_3)$
is isomorphic to that corresponding to $(\lmb_2,\lmb_1,\pm\lmb_3)$
(the $\pm$ sign depending on $\chi_2(h_1)$).
\end{rem}

\begin{rem}\label{rm:rnk}
In some cases $M_1=M_2$, $h_1^2\neq 1$ and $\chi_1^2=1$.
In those cases, let $W$ be the linear span of $a_1,a_2$.
We can change the basis $\{a_1,a_2\}$ by any other basis of $W$, and
(similarly to the case $W_3^1$) we have a quadratic form
$f:W\to\k(h_1^2-1),\ f(a)=a^2$. As in that case, $A$ is characterized by
the rank of this form, whence there shall be $3$ liftings.
\end{rem}

\begin{rem}\label{rm:cin}
In the case $V_5^4$ (and $V_6^2,\ V_6^4$ for $\Gm=C_8$) there is
no restriction to $\lmb_1$, $\lmb_2$ and
$\lmb_3$. If we replace $\{a_1,a_2\}$ by $\{sa_1,ta_2\}$
$(s,t\in\unid\k)$, then we get an isomorphic algebra but with
$(\lmb_1,\lmb_2,\lmb_3)$ replaced by $(s^2\lmb_1,t^2\lmb_2,st\lmb_3)$.
Also, \ref{rm:itr} applies here. If $\lmb_1\lmb_2\neq 0$, taking
$s^2=\lmb_2/\lmb_1$, $t=s^{-1}$, we already get that the algebra
corresponding to $(\lmb_1,\lmb_2,\lmb_3)$ is isomorphic to the one
corresponding to $(\lmb_2,\lmb_1,\lmb_3)$. For $\lmb_1\lmb_2=0$, though,
we get new relations.
There are thus infinitely many liftings, which are parametrized
by the quotient variety $\k^3/\sim$, where
\begin{align*}
&(\lmb_1,\lmb_2,\lmb_3)\sim (s^2\lmb_1,t^2\lmb_2,st\lmb_3),
	\qquad s,t\in\unid\k. \\
&(1,0,1)\sim (0,1,1) \\
&(1,0,0)\sim (0,1,0)
\end{align*}
\end{rem}

\begin{rem}\label{rm:est}
If for some reason (for instance \ref{rm:ebu} or \eqref{eq:cd3})
$\lmb_1=0$ or $\lmb_2=0$ then taking suitable scalar multiples
of $a_1$ or $a_2$ we can suppose also that $\lmb_3\in\{0,1\}$.
\end{rem}

Here is a list of the liftings up to isomorphism. We put
$$a_3=a_1a_2-\chi_2(h_1)a_2a_1.$$
$$\begin{array}{|c|c|c|c|c|c|l|} \hline
\multicolumn{7}{|c|}{\mbox{TABLE 11:
	Rank $2$, $\Gm=C_2\times C_4$, liftings}} \\ \hline
V     & \mbox{Lift.} & -b_{12} & a_1^2 & a_2^2 &
	a_3 & \mbox{Remarks} \\ \hline
V_5^1 & 3 & + & 0 & 0 & 0 & \mbox{See \ref{rm:rnk}} \\ \cline{4-6}
&&& g_2^2-1 & 0 & 0 & \\ \cline{4-6}
&&& g_2^2-1 & g_2^2-1 & 0 & \\ \hline
V_5^2 & 4 & + & 0 & 0 & 0 & \mbox{See \ref{eq:cd4}} \\ \cline{4-6}
&&& g_2^2-1 & 0 & 0 & \\ \cline{4-6}
&&& 0 & g_2^2-1 & 0 & \\ \cline{4-6}
&&& g_2^2-1 & g_2^2-1 & 0 & \\ \hline
V_5^3 & 3 & + &  0 & 0 & 0 & \mbox{See \ref{rm:rnk}} \\ \cline{4-6}
&&& g_2^2-1 & 0 & 0 & \\ \cline{4-6}
&&& g_2^2-1 & g_2^2-1 & 0 & \\ \hline
V_5^4 & \infty & + &
	\multicolumn{3}{c|}{\lmb_1(g_2^2-1),\lmb_2(g_2^2-1),\lmb_3(g_1g_2^2-1)}
	& \mbox{See \ref{rm:cin}} \\
	&&&\multicolumn{3}{c|}{(\lmb_1,\lmb_2,\lmb_3)/\sim} & \\ \hline
V_5^6 & 3 & - & 0 & 0 & 0 &
	\mbox{See \ref{eq:cd4} and \ref{rm:itr}} \\ \cline{4-6}
&&& g_2^2-1 & 0 & 0 & \\ \cline{4-6}
&&& g_2^2-1 & g_2^2-1 & 0 & \\ \hline
V_5^7 & 3 & + & 0 & 0 & 0 &
	\mbox{See \ref{rm:ebu} and \ref{rm:itr}} \\ \cline{4-6}
&&& g_2^2-1 & 0 & 0 & \\ \cline{4-6}
&&& g_2^2-1 & g_2^2-1 & 0 & \\ \hline
V_5^8 & 4 & + & 0 & 0 & 0 & \mbox{See \ref{rm:ebu}} \\ \cline{4-6}
&&& g_2^2-1 & 0 & 0 & \\ \cline{4-6}
&&& 0 & g_2^2-1 & 0 & \\ \cline{4-6}
&&& g_2^2-1 & g_2^2-1 & 0 & \\ \hline
V_5^{9} & 3 & + & 0 & 0 & 0 & \mbox{See \ref{rm:itr}} \\ \cline{4-6}
&&& g_2^2-1 & 0 & 0 & \\ \cline{4-6}
&&& g_2^2-1 & g_2^2-1 & 0 & \\ \hline
V_5^{10} & 2 & - & 0 & 0 & 0 & 
	\mbox{See \ref{rm:ebu} and \ref{eq:cd4}} \\ \cline{4-6}
&&& g_2^2-1 & 0 & 0 & \\ \hline
\multicolumn{7}{|c|}{\mbox{Continues in next page}} \\ \hline
\end{array}$$
$$\begin{array}{|c|c|c|c|c|c|l|}
\hline
V        & \mbox{Lift.} & -b_{12} & a_1^2 & a_2^2 &
	a_3 & \mbox{Remarks} \\ \hline
V_5^{13} & 4 & + & 0 & 0 & 0 & \mbox{See \ref{rm:ebu} and \ref{rm:est}} \\ \cline{4-6}
&&& g_2^2-1 & 0 & 0 & \\ \cline{4-6}
&&& 0 & 0 & g_1g_2-1 & \\ \cline{4-6}
&&& g_2^2-1 & 0 & g_1g_2-1 & \\ \hline
V_5^{22} & 1 & + & 0 & 0 & 0 & \mbox{See \ref{rm:ebu}} \\ \hline
V_5^{23} & 1 & + & 0 & 0 & 0 & \mbox{See \ref{rm:ebu}} \\ \hline
V_5^{24} & 1 & + & 0 & 0 & 0 & \mbox{See \ref{rm:ebu}} \\ \hline
V_5^{25} & 1 & - & 0 & 0 & 0 & 
	\mbox{See \ref{rm:ebu} and \ref{eq:cd4}} \\ \hline
V_5^{27} & 1 & + & 0 & 0 & 0 &
	\mbox{See \ref{rm:ebu} and \ref{eq:cd4}} \\ \hline
V_5^{28} & 2 & + & 0 & 0 & 0 & \mbox{See \ref{rm:ebu} and \ref{rm:est}} \\ \cline{4-6}
&&& 0 & 0 & g_1g_2^2-1 & \\ \hline
V_5^{30} & 1 & - & 0 & 0 & 0 &
	\mbox{See \ref{rm:ebu} and \ref{eq:cd4}} \\ \hline
V_5^{32} & 1 & + & 0 & 0 & 0 & \mbox{See \ref{rm:ebu}} \\ \hline
V_5^{33} & 1 & + & 0 & 0 & 0 & \mbox{See \ref{rm:ebu}} \\ \hline
V_5^{34} & 1 & + & 0 & 0 & 0 & \mbox{See \ref{rm:ebu}} \\ \hline
V_5^{35} & 1 & + & 0 & 0 & 0 & \mbox{See \ref{rm:ebu}} \\ \hline
V_5^{36} & 1 & + & 0 & 0 & 0 & \mbox{See \ref{rm:ebu}} \\ \hline
V_5^{37} & 1 & - & 0 & 0 & 0 &
	\mbox{See \ref{rm:ebu} and \ref{eq:cd4}} \\ \hline
V_5^{39} & 1 & + & 0 & 0 & 0 & \mbox{See \ref{rm:ebu} and \ref{rm:est}} \\ \cline{4-6}
&&& 0 & 0 & g_2^2-1 & \\ \hline
V_5^{40} & 1 & + & 0 & 0 & 0 &
	\mbox{See \ref{rm:ebu} and \ref{eq:cd4}} \\ \hline
\end{array}
$$

\subsection{$\Gm=C_8$} \label{sz8}
We have $\Aut(\Gm)\simeq C_2\times C_2$, generated by
$$f_1=(g\mapsto g^3),\quad f_2=(g\mapsto g^5).$$
The orbits of $\Gm$ are $\{g,g^3,g^5,g^7\}$, $\{g^2,g^6\}$ and $\{g^4\}$.
There are $12$ orbits of irreducible \ydm s giving Nichols algebras of
dimension a power of $2$. We choose an element in each:
$$\begin{array}{|c|l|l|c|} \hline
\multicolumn{4}{|c|}{\mbox{TABLE 12:
	modules$/\Aut(\Gm)$, $\Gm=C_8$}} \\ \hline
M(h,\chi) & h & \chi & \dim\toba M \\ \hline
Y_6^1 & g & \hat g              & 8 \\ \hline
Y_6^2 & g & \hat g^3            & 8 \\ \hline
Y_6^3 & g & \hat g^5            & 8 \\ \hline
Y_6^4 & g & \hat g^7            & 8 \\ \hline
Y_6^5 & g & \hat g^2            & 4 \\ \hline
Y_6^6 & g & \hat g^6            & 4 \\ \hline
\multicolumn{4}{|c|}{\mbox{Continues in next page}} \\ \hline
\end{array}$$
$$\begin{array}{|c|l|l|c|} \hline
M(h,\chi) & h & \chi & \dim\toba M \\ \hline
Y_6^7 & g & \hat g^4            & 2 \\ \hline
Y_6^8 & g^2 & \hat g            & 4 \\ \hline
Y_6^9 & g^2 & \hat g^3          & 4 \\ \hline
Y_6^{10} & g^2 & \hat g^2       & 2 \\ \hline
Y_6^{11} & g^2 & \hat g^6       & 2 \\ \hline
Y_6^{12} & g^4 & \hat g         & 2 \\ \hline
\end{array}$$
Hence there are $4$ different $4$-dimensional Nichols algebras
of rank $1$.
Now, for rank $2$, we suppose $M_1$ is one of $Y_6^7$, $Y_6^{10}$,
$Y_6^{11}$ or $Y_6^{12}$. We have the following possibilities (we
take into account \ref{lm:iso}).
$$\begin{array}{|c|l|l|l|l|r|r|r|r|l|c|} \hline
\multicolumn{11}{|c|}{\mbox{TABLE 13:
	Rank $2$, $\Gm=C_8$}} \\ \hline
V       & h_1 & h_2 & \chi_1 & \chi_2 & b_{11} 
	& b_{12} & b_{21} & b_{22} & b_{ij} & \dim\toba V \\ \hline
V_6^1 & g   & g   & \hat g^4 & \hat g^4         & -1 & -1 & -1 & -1
	& b^1_{-1} & 4 \\ \hline
V_6^2 & g   & g^3 & \hat g^4 & \hat g^4         & -1 & -1 & -1 & -1
	& b^1_{-1} & 4 \\ \hline
V_6^3 & g   & g^7 & \hat g^4 & \hat g^4         & -1 & -1 & -1 & -1
	& b^1_{-1} & 4 \\ \hline
V_6^4 & g   & g^5 & \hat g^4 & \hat g^4         & -1 & -1 & -1 & -1
	& b^1_{-1} & 4 \\ \hline
V_6^5 & g   & g^2 & \hat g^4 & \hat g^2         & -1 & \rmu & 1 & -1
	& [b^3_{++}]^{\rho} & 16 \\ \hline
V_6^6 & g   & g^2 & \hat g^4 & \hat g^6         & -1 & -\rmu & 1 & -1
	& [b^3_{-+}]^{\rho} & 16 \\ \hline
V_6^7 & g   & g^6 & \hat g^4 & \hat g^2         & -1 & \rmu & 1 & -1
	& [b^3_{++}]^{\rho} & 16 \\ \hline
V_6^8 & g   & g^6 & \hat g^4 & \hat g^6         & -1 & -\rmu & 1 & -1
	& [b^3_{-+}]^{\rho} & 16 \\ \hline
V_6^9    & g   & g^4 & \hat g^4 & \hat g        & -1 & \xi   & 1 & -1
	& [b^6_1]^{\rho} & 32 \\ \hline
V_6^{10} & g   & g^4 & \hat g^4 & \hat g^3      & -1 & \xi^3 & 1 & -1
	& [b^6_3]^{\rho} & 32 \\ \hline
V_6^{11} & g   & g^4 & \hat g^4 & \hat g^7      & -1 & \xi^7 & 1 & -1
	& [b^6_7]^{\rho} & 32 \\ \hline
V_6^{12} & g   & g^4 & \hat g^4 & \hat g^5      & -1 & \xi^5 & 1 & -1
	& [b^6_5]^{\rho} & 32 \\ \hline
V_6^{13} & g^2 & g^2 & \hat g^2 & \hat g^2      & -1 & -1 & -1 & -1
	& b^1_{-1} & 4 \\ \hline
V_6^{14} & g^2 & g^2 & \hat g^2 & \hat g^6      & -1 & -1 & -1 & -1
	& b^1_{-1} & 4 \\ \hline
V_6^{15} & g^2 & g^2 & \hat g^6 & \hat g^6      & -1 & -1 & -1 & -1
	& b^1_{-1} & 4 \\ \hline
V_6^{16} & g^2 & g^6 & \hat g^2 & \hat g^6      & -1 & -1 & -1 & -1
	& b^1_{-1} & 4 \\ \hline
V_6^{17} & g^2 & g^6 & \hat g^6 & \hat g^6      & -1 & -1 & -1 & -1
	& b^1_{-1} & 4 \\ \hline
V_6^{18} & g^2 & g^6 & \hat g^6 & \hat g^2      & -1 & -1 & -1 & -1
	& b^1_{-1} & 4 \\ \hline
V_6^{19} & g^4 & g^2 & \hat g   & \hat g^2      & -1 &  1 & \rmu & -1
	& b^3_{++} & 16 \\ \hline
V_6^{20} & g^4 & g^2 & \hat g   & \hat g^6      & -1 &  1 & \rmu & -1
	& b^3_{++} & 16 \\ \hline
V_6^{21} & g^4 & g^6 & \hat g   & \hat g^6      & -1 &  1 & -\rmu & -1
	& b^3_{-+} & 16 \\ \hline
V_6^{22} & g^4 & g^6 & \hat g   & \hat g^2      & -1 &  1 & -\rmu & -1
	& b^3_{-+} & 16 \\ \hline
\multicolumn{11}{|c|}{\mbox{Continues in next page}} \\ \hline
\end{array}$$
$$\begin{array}{|c|l|l|l|l|r|r|r|r|l|c|} \hline
V       & h_1 & h_2 & \chi_1 & \chi_2 & b_{11} 
	& b_{12} & b_{21} & b_{22} & b_{ij} & \dim\toba V \\ \hline
V_6^{23} & g^4 & g^4 & \hat g   & \hat g        & -1 & -1 & -1 & -1
	& b^1_{-1} & 4 \\ \hline
V_6^{24} & g^4 & g^4 & \hat g   & \hat g^3      & -1 & -1 & -1 & -1
	& b^1_{-1} & 4 \\ \hline
V_6^{25} & g^4 & g^4 & \hat g   & \hat g^7      & -1 & -1 & -1 & -1
	& b^1_{-1} & 4 \\ \hline
V_6^{26} & g^4 & g^4 & \hat g   & \hat g^5      & -1 & -1 & -1 & -1
	& b^1_{-1} & 4 \\ \hline
\end{array}$$

We compute now the liftings. For rank $1$, we have the following
algebras, where we take $a$ as in \ref{al:uodl}.
$$\begin{array}{|c|l|l|c|l|} \hline
\multicolumn{5}{|c|}{\mbox{TABLE 14:
	Rank $1$, $\Gm=C_8$}} \\ \hline
V & h & \chi & a^4 &  \mbox{Remarks} \\ \hline
Y_6^5 & g & \hat g^2 & 0 & \\ \cline{4-4}
&&& g^4-1 & \\ \hline
Y_6^6 & g & \hat g^6 & 0 & \\ \cline{4-4}
&&& g^4-1 & \\ \hline
Y_6^8 & g^2 & \hat g & 0 & \mbox{See \eqref{rm:ebu}} \\ \hline
Y_6^9 & g^2 & \hat g^3 & 0 & \mbox{See \eqref{rm:ebu}} \\ \hline
\end{array}$$
	   
\medskip
For rank $2$, as in \ref{al:uodl}, we have $a_1,a_2$ with
$a_i\in\primg A{h_i}$, $ta_it^{-1}=\chi_i(t)a_i$,
$$a_i^2=\lmb_i(h_i^2-1)\ (i=1,2),\quad a_1a_2+a_2a_1=\lmb_3(h_1h_2-1)$$
(notice that $b_{12}=\chi_2(h_1)=-1$ in all the $4$-dimensional cases).
Equations \eqref{eq:cd3}, \eqref{eq:cd4} and remarks \ref{rm:itr},
\ref{rm:rnk}, \ref{rm:cin} and \ref{rm:est} apply here also. We put
$$a_3=a_1a_2+a_2a_1.$$

$$\begin{array}{|c|c|c|c|c|l|} \hline
\multicolumn{6}{|c|}{\mbox{TABLE 15:
	Rank $2$, $\Gm=C_8$, liftings}} \\ \hline
V     & \mbox{Lift.} & a_1^2 & a_2^2 & a_3 & \mbox{Remarks} \\ \hline
V_6^1 & 3 &  0 & 0 & 0 & \mbox{See \ref{rm:rnk}} \\ \cline{3-5}
&& g^2-1 & 0 & 0 & \\ \cline{3-5}
&& g^2-1 & g^2-1 & 0 & \\ \hline
V_6^2 & \infty & \multicolumn{3}{c|}%
	{\lmb_1(g^2-1),\lmb_2(g^6-1),\lmb_3(g^4-1)}
	& \mbox{See \ref{rm:cin}} \\
&& \multicolumn{3}{c|}{(\lmb_1,\lmb_2,\lmb_3)/\sim} & \\ \hline
V_6^3 & 3 &  0 & 0 & 0 & \mbox{See \ref{rm:itr}} \\ \cline{3-5}
&& g^2-1 & 0 & 0 & \\ \cline{3-5}
&& g^2-1 & g^6-1 & 0 & \\ \hline
V_6^4 & \infty & \multicolumn{3}{c|}%
	{\lmb_1(g^2-1),\lmb_2(g^2-1),\lmb_3(g^6-1)}
	& \mbox{See \ref{rm:cin}} \\
&& \multicolumn{3}{c|}{(\lmb_1,\lmb_2,\lmb_3)/\sim} & \\ \hline
\multicolumn{6}{|c|}{\mbox{Continues in next page}} \\ \hline
\end{array}$$
$$\begin{array}{|c|c|c|c|c|l|} \hline
V     & \mbox{Lift.} & a_1^2 & a_2^2 & a_3 & \mbox{Remarks} \\ \hline
V_6^{13} & 1 &  0 & 0 & 0 &
	\mbox{See \ref{eq:cd3} and \ref{eq:cd4}}\\ \hline
V_6^{14} & 2 &  0 & 0 & 0 & \mbox{See \ref{eq:cd3} and \ref{rm:est}}\\ \cline{3-5}
&& 0 & 0 & g^4-1 & \\ \hline
V_6^{15} & 1 &  0 & 0 & 0 &
	\mbox{See \ref{eq:cd3} and \ref{eq:cd4}}\\ \hline
V_6^{16} & 1 &  0 & 0 & 0 &
	\mbox{See \ref{rm:ebu} and \ref{eq:cd3}}\\ \hline
V_6^{17} & 1 &  0 & 0 & 0 &
	\mbox{See \ref{rm:ebu} and \ref{eq:cd3}}\\ \hline
V_6^{18} & 1 &  0 & 0 & 0 &
	\mbox{See \ref{eq:cd3} and \ref{eq:cd4}}\\ \hline
V_6^{23} & 1 &  0 & 0 & 0 &
	\mbox{See \ref{rm:ebu}}\\ \hline
V_6^{24} & 1 &  0 & 0 & 0 &
	\mbox{See \ref{rm:ebu}}\\ \hline
V_6^{25} & 1 &  0 & 0 & 0 &
	\mbox{See \ref{rm:ebu}}\\ \hline
V_6^{26} & 1 &  0 & 0 & 0 &
	\mbox{See \ref{rm:ebu}}\\ \hline
\end{array}$$

\bigskip
\subsection{$\Gm=\DD_4$} \label{sd4}
We take $\DD_4$ generated by $r$ and $\sgm$ with relations
$r^4=\sgm^2=1$, $\sgm r=r^3\sgm$.
The conjugacy classes of $\Gm$ are $\{1\}$, $\{r^2\}$,
$\{\sgm,r^2\sgm\}$, $\{r\sgm,r^3\sgm\}$ and $\{r,r^3\}$.
As explained in \ref{al:casos}, there are two possibilities for $V$.
For case \ref{al:casos}\eqref{ca:d22},
$h$ must be in a conjugacy class with $2$ elements,
and hence we can take $h=\sgm$, $h=r\sgm$ or $h=r$.
However, the cases $h=\sgm$ and $h=r\sgm$ give isomorphic algebras when
bosonized by means of the morphism $(r\mapsto r,\ \sgm\mapsto r\sgm)$, and
hence we take only the cases $h=\sgm$ and $h=r$.
The commutator subgroups are
$$\Gm_{\sgm}=\{1,\sgm,r^2,r^2\sgm\}\simeq C_2\times C_2,\qquad
\Gm_{r}=\{1,r,r^2,r^3\}\simeq C_4.$$

Case \ref{al:casos}\eqref{ca:d23} arises taking $h=r^2$ and
$\rho=\rho_0$ the irreducible representation given by
$$\rho_0(r)=\mdpd 0{-1}10,\quad \rho_0(\sgm)=\mdpd 0110.$$
The following is the list of irreducible \ydm s giving
Nichols algebras of dimension $4$:
$$\begin{array}{|c|c|cc|c|c|} \hline
\multicolumn{6}{|c|}{\mbox{TABLE 16:
	$\Gm=\DD_4$}} \\ \hline
M(h,\rho) & h & \multicolumn{2}{|c|}{\rho} & b_{12}=b_{21}
	& \dim\toba M \\ \hline
Y_7^1 & \sgm & \sgm\mapsto -1,\ & r^2\mapsto 1 & -1 & 4 \\ \hline
Y_7^2 & \sgm & \sgm\mapsto -1,\ & r^2\mapsto -1 & \fsm 1 & 4 \\ \hline
Y_7^3 & r & \multicolumn{2}{|c|}{r\mapsto -1} & -1 & 4 \\ \hline
Y_7^4 & r^2 & \multicolumn{2}{|c|}{\rho_0} & -1 & 4 \\ \hline
\end{array}$$

\bigskip
We compute the liftings. Though we can state some general conditions
as in \eqref{eq:cd4}, it turns out to be simpler to use the diamond
lemma in each case.

\begin{algo} $V=Y_7^1$.
As in \ref{ta:nal1}, we have $a\in\primg A{\sgm}$ and
$b=rar^3\in\primg A{r^2\sgm}$. Since
\begin{align*}
\sgm\ganda x&=-x,\quad r\ganda x=y,\quad r\ganda y=r^2\ganda x=x, \\
\sgm\ganda y&=\sgm r\ganda x=r^3\sgm\ganda x=r(r^2\sgm)\ganda x
	=-r\ganda x=-y,\ 
\end{align*}
we have
$$\sgm a\sgm=-a,\quad rar^3=b,\quad
	\sgm b\sgm=-b,\quad rbr^3=a.$$
By \ref{rm:ebu}, $a^2=b^2=0$, and $ab+ba\in\primg A{r^2}$, from where
(taking a suitable scalar multiple of $a$ or $b$), we have
$ab+ba=\lmb(r^2-1)$, with $\lmb\in\{0,1\}$. There are hence two liftings:
\begin{enumerate}
\item $a^2=b^2=0$, $ab+ba=0$. This is the bosonization $\toba{Y_7^1}\#\k\Gm$.
\item $a^2=b^2=0$, $ab+ba=(r^2-1)$.
\end{enumerate}
\end{algo}

\begin{algo} $V=Y_7^2$.
It is similar to the case $Y_7^1$: we have $a\in\primg A{\sgm}$ and 
$b=rar^3\in\primg A{r^2\sgm}$ such that
$$\sgm a\sgm=-a,\quad rar^3=b,\quad
\sgm b\sgm=b,\quad rbr^3=-a.$$
As in that case, $a^2=b^2=0$ by \ref{rm:ebu}, and $ab-ba=\lmb(r^2-1)$
with $\lmb\in\{0,1\}$, but now the diamond lemma applies:
$$0=bba=bab-\lmb b(r^2-1)=abb-\lmb(r^2-1)b-\lmb b(r^2-1)=2\lmb b,$$
from where $\lmb=0$ for $A$ to be $32$-dimensional.

There is hence only one lifting: the bosonization $\toba{Y_7^2}\#\k\Gm$.
\end{algo}

\begin{algo} $V=Y_7^3$.
By \ref{ta:nal1}, we have here $a\in\primg Ar$,
$b=\sgm a\sgm\in\primg A{r^3}$, and
$$rar^3=-a,\quad rbr^3=-b,\quad
	\sgm a\sgm=b,\quad \sgm b\sgm=a.$$
We have now $a^2=\lmb_1(r^2-1)$, $b^2=\lmb_2(r^2-1)$. Since
$ab+ba\in\prim A$ then $ab+ba=0$ by \ref{rm:ebu}.
We can take $\lmb_i\in\{0,1\}$, from where we have four
liftings. However, the automorphism $r\mapsto r^3,\ \sgm\mapsto\sgm$
gives an isomorphism between the algebra corresponding to
$(\lmb_1=0,\lmb_2=1)$ and that of $(\lmb_1=1,\lmb_2=0)$, and thus there
are three non isomorphic liftings:
\begin{enumerate}
\item $a^2=0$, $b^2=0$, $ab+ba=0$. This is the bosonization
	$\toba{Y_7^3}\#\k\Gm$.
\item $a^2=(r^2-1)$, $b^2=0$, $ab+ba=0$.
\item $a^2=(r^2-1)$, $b^2=(r^2-1)$, $ab+ba=0$.
\end{enumerate}
\end{algo}

\begin{algo} $V=Y_7^4$.
As in \ref{ta:nal2} (with $t=r$), we have here $a,b\in\primg A{r^2}$, with
$$rar^3=b,\quad rbr^3=-a,\quad
	\sgm a\sgm=b,\quad \sgm b\sgm=a.$$
By \ref{rm:ebu}, we have $a^2=b^2=ab+ba=0$, and thus
there is only one lifting: the bosonization $\toba{Y_7^4}\#\k\Gm$.
\end{algo}

\subsection{$\Gm=\HH$} \label{sh}
We take $\HH=\{e,(-e),\ii,-\ii,\jj,-\jj,\kk,-\kk\}$ with the
standard relations $\ii\jj=\kk$, etc. The conjugacy
classes are $\{e\}$, $\{(-e)\}$, $\{\ii,-\ii\}$, $\{\jj,-\jj\}$ and
$\{\kk,-\kk\}$. Since the last three classes give isomorphic algebras via
the morphisms
$$(\ii\mapsto \jj\mapsto \kk\mapsto \ii)\quad \mbox{and}\quad
	(\ii\mapsto -\ii,\ \jj\mapsto -\jj,\ \kk\mapsto -\kk),$$
there is only one case as in \ref{al:casos}\eqref{ca:d22}, which
is given by $h=\ii$, $\chi(\ii)=-1$.

\noindent There is also only one case as in \ref{al:casos}\eqref{ca:d23},
which is given by $h=(-e)$, $\rho=\rho_0$, where
$$\rho_0(\ii)=\mdpd 0{-1}10,\ \rho_0(\jj)=\mdpd{\rmu}00{-\rmu},\ 
	\rho_0(\kk)=\mdpd 0{\rmu}{\rmu}0.$$

The list of $4$-dimensional Nichols algebras giving non isomorphic algebras
when bosonized is given hence by
$$\begin{array}{|c|c|c|c|c|} \hline
\multicolumn{5}{|c|}{\mbox{TABLE 17:
	$\Gm=\HH$}} \\ \hline
M(h,\rho) & h & \rho & b_{12}=b_{21} & \dim\toba M \\ \hline
Y_8^1 & \ii & \ii\mapsto -1 & -1 & 4 \\ \hline
Y_8^2 & (-e) & \rho_0 & -1 & 4 \\ \hline
\end{array}$$

The liftings are:
\begin{algo} $V=Y_8^1$.
As in \ref{ta:nal1} we have $a\in\primg A{\ii}$ and
$b=\jj a(-\jj)\in\primg A{-\ii}$, with the action given by
\begin{alignat*}{3}
\ii a(-\ii) &=-a,\quad &&\jj a(-\jj)=b,\quad &&\kk a(-\kk)=-b, \\
\ii b(-\ii) &=-b,\quad &&\jj b(-\jj)=a,\quad &&\kk b(-\kk)=-a.
\end{alignat*}
We have $\lmb_1,\lmb_2\in\{0,1\}$ such that $a^2=\lmb_1((-e)-1)$,
$b^2=\lmb_2((-e)-1)$, $ab+ba\in\prim A\so ab+ba=0$. The diamond lemma
says, however, that $\lmb_1=\lmb_2$, since
\begin{align*}
\jj a^2&=\lmb_1\jj((-e)-1)=\lmb_1((-\jj)-\jj), \\
\jj a^2&=b\jj a=b^2\jj=\lmb_2((-e)-1)\jj=\lmb_2((-\jj)-\jj).
\end{align*}
Hence, there are two liftings:
\begin{enumerate}
\item $a^2=0$, $b^2=0$, $ab+ba=0$. This is the bosonization
	$\toba{Y_8^1}\#\k\Gm$.
\item $a^2=((-e)-1)$, $b^2=((-e)-1)$, $ab+ba=0$.
\end{enumerate}
\end{algo}

\begin{algo} $V=Y_8^2$.
As in \ref{ta:nal2} we have $a,b\in\primg A{-e}$ with the action given by
\begin{alignat*}{3}
\ii a(-\ii)&= b,\quad &&\jj a(-\jj)= \rmu a,\quad &&\kk a(-\kk)=\rmu b, \\
\ii b(-\ii)&=-a,\quad &&\jj b(-\jj)=-\rmu b,\quad &&\kk b(-\kk)=\rmu a.
\end{alignat*}
By \ref{rm:ebu}, we have $a^2=b^2=0$, and $ab+ba=0$, and hence there is
only one lifting: the bosonization $\toba{Y_8^2}\#\k\Gm$.
\end{algo}

\section{$\Gm$ of order $16$}
Now $\dim\toba V=2$, from where we must have $V=M(h,\chi)$, $h\in Z(\Gm)$
and $\chi$ a character such that $\chi(h)=-1$. Hence, for $\Gm$ to have
\ydm s $V$ with $\dim\toba V=2$, it is necessary that
$Z(\Gm)\not\subset [\Gm;\Gm]$. There are $14$ groups of order $16$. Among
them, $3$ do not verify this condition.

As explained in \ref{ta:lev}, we take $a\in\primg Ah$ such that
$ga=\chi(g)ag$ for $g\in\Gm$. Furthermore, there exists
$\lmb\in\{0,1\}$ such that $a^2=\lmb(h^2-1)$.
Analogously to \eqref{eq:cd2}, for $A$ to be $32$-dimensional we have
$$\lmb g(h^2-1)=ga^2=\chi^2(g)a^2g=\lmb\chi^2(g)(h^2-1)g
	\quad\forall g\in\Gm$$
whence $\lmb(\chi^2(g)-1)(h^2-1)=0$. Combining this with \ref{rm:ebu}
we get
\begin{equation}\label{eq:con4}
\lmb(\chi^2-1)=0.
\end{equation}

\subsection{$\Gm=C_2\times C_2\times C_2\times C_2$} \label{sz2z2z2z2}
As in the case $\Gm=C_2\times C_2\times C_2$, the group $\Aut(\Gm)$ acts
transitively on the \ydm s $V$ such that $\dim\toba V=2$,
from where we get only one algebra after
bosonization. We take hence $V=M(g_1,\hat g_1)$. By \ref{rm:ebu},
there is only one lifting: the bosonization $\toba V\#\k\Gm$.

\subsection{$\Gm=C_2\times C_2\times C_4$} \label{sz2z2z4}
There are $4$ orbits with respect to the action of $\Aut(\Gm)$:
$\{1\}$, $\{g_3^2\}$, $\{g_1^ig_2^jg_3^k,\ (i,j)\neq(0,0),\ k=0,2\}$,
$\{g_1^ig_2^jg_3^k,\ k=1,3\}$ (we are using the notation at the end
of \ref{ss:bdac}). The first orbit, as always, does not give a
Nichols algebra of dimension $2$. For the others, we have the
following list of \ydm s giving non isomorphic algebras when
bosonized (we indclude only the cases $\dim\toba V=2$).
We include also their corresponding liftings:
$$\begin{array}{|c|c|c|c|c|l|} \hline
\multicolumn{6}{|c|}{\mbox{TABLE 18:
	$\Gm=C_2\times C_2\times C_4$, liftings}} \\ \hline
M(h,\chi) & \mbox{Liftings} & h & \chi & a^2 & \mbox{Remarks} \\ \hline
Y_9^1 & 1 & g_3^2 & \hat g_3 & 0 & \mbox{See \ref{rm:ebu}} \\ \hline
Y_9^2 & 1 & g_1 & \hat g_1 & 0 & \mbox{See \ref{rm:ebu}} \\ \hline
Y_9^3 & 1 & g_1 & \hat g_1\hat g_3 & 0 & \mbox{See \ref{rm:ebu}} \\ \hline
Y_9^4 & 2 & g_3 & \hat g_3^2 & 0 & \\ \cline{5-5}
&&&& (g_3^2-1) & \\ \hline
Y_9^5 & 2 & g_3 & \hat g_1\hat g_3^2 & 0 & \\ \cline{5-5}
&&&& (g_3^2-1) & \\ \hline
\end{array}$$

\subsection{$\Gm=C_4\times C_4$} \label{sz4z4}
The orbits are $\{1\}$, $\{g_1^2,g_1^2g_2^2,g_2^2\}$,
$\{g_1^ig_2^j,\ \ 2\not{|}i\mbox{ or }2\not{|}j\}$ (we are using the
notation at the end of \ref{ss:bdac}). The possibilities are
$$\begin{array}{|c|c|c|c|c|l|} \hline
\multicolumn{6}{|c|}{\mbox{TABLE 19:
	$\Gm=C_4\times C_4$, liftings}} \\ \hline
M(h,\chi) & \mbox{Liftings} & h & \chi & a^2 & \mbox{Remarks} \\ \hline
Y_{10}^1 & 2 & g_1 & \hat g_1^2 & 0 & \\ \cline{5-5}
&&&& (g_1^2-1) & \\ \hline
Y_{10}^2 & 1 & g_1 & \hat g_1^2\hat g_2 & 0 
		& \mbox{See \ref{eq:con4}} \\ \hline
Y_{10}^3 & 1 & g_1^2 & \hat g_1 & 0
		& \mbox{See \ref{rm:ebu}} \\ \hline
\end{array}$$

\subsection{$\Gm=C_2\times C_8$} \label{sz2z8}
The orbits are $\{1\}$, $\{g_2^4\}$, $\{g_1g_2^2,g_1g_2^6\}$,
$\{g_1,g_1g_2^4\}$, $\{g_2^2,g_2^6\}$, $\{g_1^ig_2^j,\ j=1,3,5,7\}$.
The list is
$$\begin{array}{|c|c|c|c|c|l|} \hline
\multicolumn{6}{|c|}{\mbox{TABLE 20:
	$\Gm=C_2\times C_8$, liftings}} \\ \hline
M(h,\chi) & \mbox{Liftings} & h & \chi & a^2 & \mbox{Remarks} \\ \hline
Y_{11}^1 & 1 & g_2^4 & \hat g_2 & 0 & \mbox{See \ref{rm:ebu}} \\ \hline
Y_{11}^2 & 2 & g_1g_2^2 & \hat g_1 & 0 & \\ \cline{5-5}
&&&& (g_2^4-1) & \\ \hline
Y_{11}^3 & 1 & g_1g_2^2 & \hat g_2^2 & 0
		& \mbox{See \ref{eq:con4}} \\ \hline
Y_{11}^4 & 1 & g_1 & \hat g_1 & 0 & \mbox{See \ref{rm:ebu}} \\ \hline
Y_{11}^5 & 1 & g_1 & \hat g_1\hat g_2^2 & 0 
		& \mbox{See \ref{rm:ebu}} \\ \hline
Y_{11}^6 & 1 & g_1 & \hat g_1\hat g_2 & 0 
		& \mbox{See \ref{rm:ebu}} \\ \hline
Y_{11}^7 & 1 & g_2^2 & \hat g_2^2 & 0
		& \mbox{See \ref{eq:con4}} \\ \hline
Y_{11}^8 & 1 & g_2^2 & \hat g_2^6 & 0
		& \mbox{See \ref{eq:con4}} \\ \hline
Y_{11}^9 & 1 & g_2^2 & \hat g_1\hat g_2^2 & 0
		& \mbox{See \ref{eq:con4}} \\ \hline
Y_{11}^{10} & 2 & g_2 & \hat g_2^4 & 0 & \\ \cline{5-5}
&&&& (g_2^2-1) & \\ \hline
Y_{11}^{11} & 2 & g_2 & \hat g_1\hat g_2^4 & 0 & \\ \cline{5-5}
&&&& (g_2^2-1) & \\ \hline
\end{array}$$

\subsection{$\Gm=C_{16}$} \label{sz16}
The orbits are
$$\{1\},\ \{g^8\},\ \{g^4,g^{12}\},\ \{g^2,g^6,g^{10},g^{14}\},
\ \{g^i,\ \ 2\not{|}i\}.$$
The list is
$$\begin{array}{|c|c|c|c|c|l|} \hline
\multicolumn{6}{|c|}{\mbox{TABLE 21:
	$\Gm=C_{16}$, liftings}} \\ \hline
M(h,\chi) & \mbox{Liftings} & h & \chi & a^2 & \mbox{Remarks} \\ \hline
Y_{12}^1 & 1 & g^8 & \hat g & 0
		& \mbox{See \ref{rm:ebu}} \\ \hline
Y_{12}^2 & 1 & g^4 & \hat g^2 & 0
		& \mbox{See \ref{eq:con4}} \\ \hline
Y_{12}^3 & 1 & g^4 & \hat g^6 & 0
		& \mbox{See \ref{eq:con4}} \\ \hline
Y_{12}^4 & 1 & g^2 & \hat g^4 & 0
		& \mbox{See \ref{eq:con4}} \\ \hline
Y_{12}^5 & 1 & g^2 & \hat g^{12} & 0
		& \mbox{See \ref{eq:con4}} \\ \hline
Y_{12}^6 & 2 & g & \hat g^8 & 0 & \\ \cline{5-5}
&&&& (g^2-1) & \\ \hline
\end{array}$$

\subsection{Non abelian groups} \label{sgna}
As said above, the non abelian groups which verify
$Z(\Gm)\not\subset[\Gm;\Gm]$ are $6$, and can be presented in the following
list:
$$\begin{array}{|c|c|c|c|c|} \hline
\multicolumn{5}{|c|}{\mbox{TABLE 22:
	Non abelian groups of order $16$}} \\ \hline
\mbox{Group} & \mbox{Generators, relations} & Z(\Gm) & [\Gm;\Gm] & \Gm_{ab}\\ \hline
B_1 & (g_1,g_2) / (g_1^8,g_2^2,g_2g_1g_2g_1^3) & (g_1^2) 
	& (g_1^4) & C_4\times C_2 \\ \hline
B_2 & (g_1,g_2,g_3) / R_2 & (g_1)
	& (g_1^2) & C_2\times C_2\times C_2 \\ \hline
B_3 & (g_1,g_2) / (g_1^4,g_2^4,g_2g_1g_2^3g_1) & (g_1^2,g_2^2) 
	& (g_1^2) & C_2\times C_4 \\ \hline
B_4=\DD_4\times C_2 & (r,\sgm,t) & (t,r^2)
	& (r^2) &  C_2\times C_2\times C_2 \\ \hline
B_5 & (g_1,g_2,g_3) / R_5 & (g_1^2,g_2) 
	& (g_2) & C_4\times C_2 \\ \hline
B_6=\HH\times C_2 & (\ii,\jj,t) & (t,-e)
	& (-e) & C_2\times C_2\times C_2 \\ \hline
\end{array}$$
where we put for $B_2$ and $B_5$ the subgroups $R_2$ and $R_5$ respectively,
generated by the elements
\begin{align*}
R_2:\quad &g_1^4,\ g_2^2,\ g_3^2,\ [g_1;g_2],\ [g_1;g_3],\ [g_3;g_2]g_1^2, \\
R_5:\quad &g_1^4,\ g_2^2,\ g_3^2,\ [g_1;g_2],\ [g_3;g_2],\ [g_1;g_3]g_2.
\end{align*}
The name $B_i$ corresponds to the classification in \cite[Sect. 118]{bur}.
In these groups, we have $[\Gm;\Gm]\subset Z(\Gm)$,
$(Z(\Gm)-[\Gm;\Gm])=\{h_1,h_2\}$. Hence, there are only two possibilities
for $h$ in each group. 

For $\Gm=B_4$, for instance, $h_1=t$ and $h_2=r^2t$.
Since the automorphism
$$(r\mapsto r,\ \sgm\mapsto\sgm,\ t\mapsto r^2t)$$
carries $h_1$ to $h_2$, we consider only $h=t$. There are $4$ characters
$\chi\in\hat\Gm$ such that $\chi(t)=-1$. They are given by
$$\chi(t)=-1,\quad \chi(r)=\pm 1,\quad \chi(\sgm)=\pm 1.$$
However, the automorphisms
$$(t\mapsto t,\ r\mapsto t^ir,\ \sgm\mapsto t^j\sgm)\quad i,j=0,1$$
act transitively on these representations,
whence we have only one algebra.

Doing the same analysis with the other groups, we arrive to the
following list (we take $a$ as explained at the beginning of the section):

$$\begin{array}{|c|c|c|c|ll|c|l|} \hline
\multicolumn{8}{|c|}{\mbox{TABLE 23:
	$\Gm$ non abelian of order $16$, liftings}} \\ \hline
\mbox{Group} & M(h,\chi) & \mbox{Lift.} & h & \multicolumn{2}{c|}{\chi}
	& a^2 & \mbox{Remarks} \\ \hline
B_1 & Y_{13}^1 & 1 & g_1^2 &
	\chi(g_1)=\rmu & \chi(g_2)=1
	& 0 & \mbox{See \ref{eq:con4}} \\ \hline
B_1 & Y_{13}^2 & 1 & g_1^2 &
	\chi(g_1)=\rmu & \chi(g_2)=-1
	& 0 & \mbox{See \ref{eq:con4}} \\ \hline
B_1 & Y_{13}^3 & 1 & g_1^2 &
	\chi(g_1)=-\rmu & \chi(g_2)=1
	& 0 & \mbox{See \ref{eq:con4}} \\ \hline
B_1 & Y_{13}^4 & 1 & g_1^2 &
	\chi(g_1)=-\rmu & \chi(g_2)=-1
	& 0 & \mbox{See \ref{eq:con4}} \\ \hline
B_2 & Y_{13}^5 & 2 & g_1 &
	\chi(g_1)=-1 & \chi(g_2)=1
	& 0 & \\ \cline{7-7}
&&&&& \chi(g_3)=1 & g_1^2-1 & \\ \hline
B_2 & Y_{13}^6 & 2 & g_1 &
	\chi(g_1)=-1 & \chi(g_2)=-1
	& 0 & \\ \cline{7-7}
&&&&& \chi(g_3)=-1 & g_1^2-1 & \\ \hline
B_3 & Y_{13}^7 & 1 & g_2^2 &
	\chi(g_1)=1 & \chi(g_2)=\rmu
	& 0 & \mbox{See \ref{rm:ebu}} \\ \hline
B_3 & Y_{13}^8 & 1 & g_1^2g_2^2 &
	\chi(g_1)=1 & \chi(g_2)=\rmu
	& 0 & \mbox{See \ref{rm:ebu}} \\ \hline
B_4 & Y_{13}^9 & 1 & t &
	\chi(t)=-1 & \chi(r)=1
	& 0 & \mbox{See \ref{rm:ebu}} \\
&&&&& \chi(\sgm)=1 && \\ \hline
B_5 & Y_{13}^{10} & 1 & g_1^2 &
	\chi(g_1)=\rmu & \chi(g_2)=1
	& 0 & \mbox{See \ref{rm:ebu}} \\
&&&&& \chi(g_3)=1 && \\ \hline
B_6 & Y_{13}^{11} & 1 & t &
	\chi(t)=-1 & \chi(\ii)=1
	& 0 & \mbox{See \ref{rm:ebu}} \\
&&&&& \chi(\ii)=1 && \\ \hline
\end{array}$$

\renewcommand{\theequation}{\arabic{section}.\arabic{equation}}
\section{The classification is complete} \label{srtoba}
We use the techniques of \cite[Sect. 8]{as2}. If $A$ is a
$32$-dimensional pointed Hopf algebra, we apply the procedures
described in \ref{al:prl} and get $R$ a braided Hopf algebra in
$\ydskg$ satisfying \eqref{toba1}, \eqref{toba2} and \eqref{toba3}.
We may consider the (braided) dual $\dad R$ as in \cite[Def. 2.2.2]{ag},
which is a braided Hopf algebra in $\ydskg$ satisfying
\eqref{toba1}, \eqref{toba2} and \eqref{toba4}. Let $W=\agg{\dad R}1$.
We get in this way the sequence of braided Hopf algebras
$$TW \to \dad R \sobre \toba W,$$
whence $\toba W\#\k\Gm$ is an algebra of dimension $2^n$,
with $n\le 5$. If $n=5$ then $\dad R=\toba W$, from where $\dad R$
satisfies \eqref{toba3} and hence $R$ satisfies
\eqref{toba4}, whence $R$ is a Nichols algebra
and $A$ is in the list above. We must prove thence that the above
surjection $\dad R\sobre\toba W$ is in fact a bijection. This is
the same as to prove that all the relations in $\toba W$ are
relations in $\dad R$, i.e., if $z\in TW$ vanishes in $\toba W$,
then it vanishes in $\dad R$. For this, notice first that for $\Gm$
non abelian there are no bosonizations of Nichols algebras
over $\Gm$ with dimension $2^n$, $n<5$ (see \ref{rm:cdna}).

Consider then $\Gm$ abelian. If $W$ is $1$-dimensional, let $x$
be a generator of $W$. We have $c(x\otimes x)=qx\otimes x$ for some
$q$ a root of unity. Let $r$ be the order of $q$ ($q\neq 1$ because
of \ref{rm:qn1}). Thus, $x^r\in\prim{\dad R}$ and
$c(x^r\otimes x^r)=(x^r\otimes x^r)$, whence by \ref{rm:qn1}
we have $x^r=0\in\dad R$. This implies that $\dad R$ verifies
\eqref{toba3}, and we are done.

\smallskip
Consider now $\dim W\ge 2$. By \ref{ex:qls}, $\dim\toba W\ge 2^{\dim W}$,
and then
$$\dim\toba W\#\k\Gm\ge 2^{\dim W}|\Gm|\ge 4|\Gm|.$$
Thus, if $R$ were not a Nichols algebra, $|\Gm|$ would be $\le 4$.

Let $\Gm=C_2$. We have $W=\oplus_{i=1}^{\te}M(g,\hat g)$.
Then for $y\in W$ we have 
$y^2\in\prim{\dad R}$ and
$c(y^2\otimes y^2)=y^2\otimes y^2$, from where by \ref{rm:qn1}
we have $y^2=0\ \forall y\in W$. These are the relations for
$\toba W$, and we are done.

If $|\Gm|=4$, the bijectivity of the projection above is a consequence
of the following
\begin{prop}
Let $\Gm$ be a finite group, $S=\oplus_i\agg Si\in\ydskg$ a graded
braided Hopf algebra that verifies \eqref{toba2} and \eqref{toba4},
and such that $\dim S=8$. Then $S$ is a Nichols algebra.
\end{prop}

\begin{pf}
Let $V=\agg S1$ and $R'=\toba V$. We have $S\sobre R'$ the canonical
projection. If $\dim R'=8$, we are done. Consider then $\dim R'<8$.
Since the projection induces an epimorphism $S\#\k\Gm\sobre R'\#\k\Gm$,
we have $\dim R'=1,\ 2$ or $4$. Since
$$\dim R\ge\dim\agg R0+\dim\agg R1=\dim\agg S0+\dim\agg S1
	=1+\dim\agg S1,$$
we have $\dim R'=2$ or $4$. If $\dim R'=2$, then $\dim\agg S1=1$ and
by the considerations for the rank-$1$ case, we would have $S=R'$.
Thus $\dim R'=4$. Now \ref{pr:dualpuan} implies that $\dim\agg {R'}1=2$,
and then by \cite[Prop 3.1.11]{ag}, $\agg{R'}1$ comes from the abelian
case. Furthermore, by the same computations as in \ref{ta:nal1}
and \ref{ta:nal2}, the matrix of $R'$ is $\mdpd{-1}{\zeta}{\zeta}{-1}$
(with $\zeta=\pm 1$) in some basis $\{x,y\}$. Now, as in \ref{ex:taft},
we have $x^2\in\prim S$ and $y^2\in\prim S$, and it is straightforward
to see that $z=xy-\zeta yx\in\prim S$. Moreover,
$$c(x^2\otimes x^2)=x^2\otimes x^2,\quad
	c(y^2\otimes y^2)=y^2\otimes y^2,\quad
	c(z\otimes z)=z\otimes z.$$
Thus, thanks to \ref{rm:qn1}, we have $x^2=y^2=z=0$. Since $S$ is
generated by $x$ and $y$, this implies $\dim S=4$, a contradiction.
\end{pf}


\begin{thebibliography}{CDRi}
\bibitem[AG]{ag} N. Andruskiewitsch \& M. Gra\~na,
	{\it Braided Hopf algebras over non abelian finite groups},
	Bolet\'\i n de la Acad. Nac. Cs. (C\'ordoba) {\bf 63} (1999), 45--78.
	Also in {\tt q-alg 9802074}.
\bibitem[AS1]{as1} N. Andruskiewitsch \& H.-J. Schneider,
	{\it Lifting of quantum linear spaces and pointed Hopf algebras
	of order $p^3$}, J.Algebra {\bf 209} (1998) 659--691.
\bibitem[AS2]{as2} N. Andruskiewitsch \& H.-J. Schneider,
	{\it Finite quantum groups and {C}artan matrices},
	to appear in Adv.Math.
\bibitem[AS3]{as3} N. Andruskiewitsch \& H.-J. Schneider,
	{\it Pointed Hopf algebras of dimension $p^4$}, in preparation.
\bibitem[Be]{be} M. Beattie
	{\it An isomorphism theorem for Ore extension Hopf algebras},
	to appear in Comm. Alg.
\bibitem[BDG]{bdg} M. Beattie, S. D$\breve{\mbox{a}}$sc$\breve{\mbox{a}}$lescu
	\& L. Gr\"unenfelder,
	{\it On the number of types of finite-dimensional Hopf algebras},
	Inventiones Math. {\bf 136} 1, (1999), 1--7.
\bibitem[Bu]{bur} W. Burnside,
	{\it Theory of groups of finite order},
	2nd edition, Dover Pub., New York, 1955.
\bibitem[CDR]{cdr} S. Caenepeel, S. D$\breve{\mbox{a}}$sc$\breve{\mbox{a}}$lescu
	\& S. Raianu,
	{\it Classifying pointed Hopf algebras of dimension $16$},
	to appear in Comm. Alg.
\bibitem[D]{dso} S. D$\breve{\mbox{a}}$sc$\breve{\mbox{a}}$lescu
	{\it Pointed Hopf algebras of dimension $p^n$ with large coradical}, preprint.
\bibitem[G]{gp5} M. Gra\~na,
	{\it On pointed Hopf algebras of dimension $p^5$}, to appear in
	Glasgow Math. Journal.
\bibitem[Ge]{ge} S. Gelaki,
	{\it On pointed Hopf algebras and Kaplansky's tenth conjecture},
	J.Algebra {\bf 209} (1998), 635--657.
\bibitem[M]{mon} S. Montgomery,
	{\it Hopf algebras and their actions on rings},
	AMS (1993), CMBS {\bf 82}.
\bibitem[N]{n} W.D. Nichols,
	{\it Bialgebras of type one},
	Comm. in Alg. {\bf 6} (1978), 1521--1552.
\bibitem[Ros]{ros} M. Rosso,
	{\it Groupes quantiques et alg{\'e}bres de battage quantiques},
	CRAS Paris {\bf 320} S\'erie I (1995), 145--148;
	{\it Quantum groups and quantum shuffles},
	Inventiones Math. {\bf 133} (1998), 399--416.
\bibitem[S]{s} M.E. Sweedler,
	{\it Hopf algebras}, Benjamin, New York, 1969.
\bibitem[Sch]{schb} P. Schauenburg,
	{\it A Characterization of the {B}orel-like Subalgebras of
	Quantum Enveloping Algebras},
	Comm. in Alg. {\bf 24} (1996), 2811--2823.
\end{thebibliography}
\end{document}